\documentclass[reqno, 12pt]{amsart}
\usepackage{amsmath}
\usepackage{amssymb}
\usepackage{amsthm}
\usepackage{mathtools}
\usepackage{tikz-cd}
\usepackage{enumerate}
\usepackage[mathscr]{eucal}
\usepackage{graphicx}

\newtheorem{theorem}{Theorem}[section]

\newtheorem{proposition}[theorem]{Proposition}
\newtheorem{corollary}[theorem]{Corollary}
\newtheorem{definition}[theorem]{Definition}
\theoremstyle{definition}
\newtheorem{example}[theorem]{Example}
\newtheorem{remark}[theorem]{Remark}

\begin{document}

	\title{Set convergences and uniform convergence of distance functionals on a bornology}
%
%
%

	\title[]{Set convergences and uniform convergence of distance functionals on a bornology}
\author{Yogesh Agarwal \and Varun Jindal}

\address{Yogesh Agarwal: Department of Mathematics, Malaviya National Institute of Technology Jaipur, Jaipur-302017, Rajasthan, India}
\email{yagarwalm247@gmail.com}

\address{Varun Jindal: Department of Mathematics, Malaviya National Institute of Technology Jaipur, Jaipur-302017, Rajasthan, India}
\email{vjindal.maths@mnit.ac.in}

\maketitle

\begin{abstract}
 For a metric space $(X,d)$, Beer, Naimpally, and Rodr{\'i}guez-L{\'o}pez in (\cite{Idealtopo}) proposed a unified approach to explore set convergences via uniform convergence of distance functionals on members of an arbitrary family $\mathcal{S}$ of subsets of $X$. The associated topology on the collection $CL(X)$ of all nonempty closed subsets of $(X,d)$ is denoted by $\tau_{\mathcal{S},d}$. As special cases, this unified approach includes classical Wijsman, Attouch-Wets, and Hausdorff distance topologies. In this article, we investigate various topological characteristics of the hyperspace $(CL(X), \tau_{\mathcal{S},d})$ when $\mathcal{S}$ is a bornology on $(X,d)$.  In order to do this, a new class of bornologies and a new metric topology on $CL(X)$ have been introduced and studied.
\end{abstract}	

\section{Introduction}
The convergence of a sequence or net of subsets of a topological space is commonly known as set convergence, and a topology on a collection of subsets of a topological space is called a hyperspace topology. There is a rich literature concerning set convergences and hyperspace topologies. Historically, most of these hyperspace topologies have been introduced by non topologists. This is precisely due to the applicability of the theory of hyperspace to various other branches of mathematics such as convex analysis, variational analysis, and optimization theory (\cite{ToCCoS, bv1996, lucchetti2006convexity, rockafellar2009variational}).

For a metric space $(X,d)$, we denote the set of all nonempty (nonempty and closed) subsets of $X$ by $\mathcal{P}_0(X)$ ($CL(X)$). The most studied hyperspace topology on $CL(X)$ is the Hausdorff distance topology (\cite{haus,pompeiu}), denoted by $\tau_{H_{d}}$,  given by the Hausdorff distance $$H_d(A,C) = \displaystyle{\sup_{x\in X}|d(x,A)-d(x,C)|} \quad (A,C \in CL(X)),$$ where $d(x, A) = \inf\{d(x,a): a \in A\}$.  For a nonempty subset $A$ of $(X,d)$, the \textit{distance functional} $d(\cdot, A)$ on $X$ corresponding to $A$ is defined by $x \to d(x, A)$ for any $x \in X$. It follows that a net $(A_{\lambda})$ in $CL(X)$ converges to $A \in CL(X)$ in the Hausdorff distance topology if and only if the net $(d(\cdot, A_{\lambda}))$ converges uniformly to $d(\cdot, A)$. The Hausdorff distance topology is useful in several framework such as fractals, image processing (\cite{hausfrac-2,hausfrac}) but not always suitable when we deal with unbounded sets.

A topology which is well-behaved even for unbounded sets and  applicable in various scenarios (\cite{AW-2, attouch1991quantitative, beer1990conjugate, rockafellar2009variational}) is the Attouch-Wets topology or bounded Hausdorff metric topology on $CL(X)$ (\cite{attouch1991topology,AW-1}), denoted by $\tau_{AW_{d}}$. The convergence corresponding to the Attouch-Wets topology  is defined as follows: a net $(A_{\lambda})$ is said to be $\tau_{AW_d}$-convergent to $A$ in $CL(X)$ provided the net $(d(\cdot, A_{\lambda}))$ converges to $d(\cdot, A)$ uniformly on bounded subsets of $(X,d)$.


In (\cite{Idealtopo}), Beer, Naimpally, and Rodr{\'i}guez-L{\'o}pez studied hyperspace convergence on $\mathcal{P}_0(X)$ defined through a unified approach by using an arbitrary family $\mathcal{S}$ of nonempty subsets of $(X,d)$. This approach incorporates convergence with respect to the Hausdorff distance topology and  Attouch-Wets topology as special cases. For a metric space $(X,d)$ and a family $\mathcal{S}$ of nonempty subsets of $X$, they defined a net $(A_{\lambda})$ of nonempty subsets of $X$ converges to $A$ in $\mathcal{P}_0(X)$ if and only if the associated net $(d(\cdot, A_{\lambda}))$ converges to $d(\cdot, A)$ uniformly on members of $\mathcal{S}$. They denoted the corresponding topology by $\tau_{\mathcal{S},d}$. As mentioned in (\cite{beer2013gap}) this unified approach was first proposed by Cornet in (\cite{cornet1973topologies}). Recently, the authors in (\cite{agarwal2024set}) have studied the relationship of upper part of $\tau_{\mathcal{S}, d}$-convergence with various other set convergences. 

In the literature, several authors have studied various topological properties such as metrizability, second countability, separability, etc., of various hyperspace topologies (\cite{caoamsterdam,wijr2,costantini2004tightness, fellFlachsmeyer, francaviglia1985quasi, michael1951}). For example, the Wijsman topology (Fell topology) on $CL(X)$ is metrizable if and only if the underlying metric space is separable (locally compact and separable) (\cite{francaviglia1985quasi, fellFlachsmeyer}). In (\cite{costantini2004tightness}), the authors studied tightness and related properties of lower Vietoris and Fell topologies. In (\cite{Idealtopo}), it has been shown that if $\mathcal{S}$ is a bornology with a countable base, then the topology  $\tau_{\mathcal{S},d}$ on $\mathcal{P}_0(X)$ is metrizable. 
 
The aim of this article is a careful investigation of the topological character of the space $(CL(X), \tau_{\mathcal{S},d})$ when $\mathcal{S}$ is a bornology on $(X,d)$. More precisely, we investigate submetrizability, metrizability, and second countability of $(CL(X), \tau_{\mathcal{S},d})$. To study these properties, we first introduce a new concept of totally bounded family on $X$ (see, Definition \ref{totally bounded family}). It is shown that the class of bornologies which form a totally bounded family on $X$ is bigger than the class of bornologies with a countable base. The relation of totally bounded family on $X$ with the concept of weakly $\mathcal{S}$-totally bounded (\cite{BCL}) is also investigated.  

The organization of the paper is as follows. In section 2, we provide basic terminologies and define the topology $\tau_{\mathcal{S},d}$. In Section 3, we explore the concept of a totally bounded family on $(X,d)$. A number of examples are given which relate this concept to various other notions related to bornologies.  In Section 4, we characterize submetrizability of the space $(CL(X), \tau_{\mathcal{S},d})$. This is done by producing explicitly a metric on $CL(X)$ whose topology is weaker than $\tau_{\mathcal{S},d}$. An independent study of $CL(X)$ equipped with this new metric is also carried out.  In Section 5, we examine the metrizability of the hyperspace $(CL(X), \tau_{\mathcal{S},d})$, and show that $(CL(X), \tau_{\mathcal{S},d})$ is metrizable if and only if $\mathcal{S}$ is a totally bounded family on $X$. The final section of this paper is devoted to investigate some countability properties of the space $(CL(X), \tau_{\mathcal{S},d})$.       

\section{Preliminaries}

For a metric space $(X,d)$, a family $\mathcal{S} \subseteq \mathcal{P}_0(X)$ which is closed under finite union, taking subsets, and forms a cover of $X$ is called a \textit{bornology} on $X$. Examples of some bornologies include: 
\begin{itemize}
	\item $\mathcal{F}(X) =$ the family of all nonempty finite subsets of $X$;
	\item  $\mathcal{B}_d(X) =$ the family of all nonempty $d$-bounded subsets of $X$;
	\item $\mathcal{K}(X) = $ the family of all nonempty relatively compact subsets of $X$;
	\item $\mathcal{T}\mathcal{B}_{d}(X) = $ the family of all nonempty $d$-totally bounded subsets of $X$.
\end{itemize}
 A family $\mathcal{S} \subseteq \mathcal{P}_0(X)$ closed under finite union and taking subsets is called an \textit{ideal} on $X$. A base of a bornology $\mathcal{S}$ is a cofinal subfamily, $\mathcal{S}_0$ of $\mathcal{S}$, that is, for each $S \in \mathcal{S}$, there is a $S_0 \in \mathcal{S}_0$ such that $S \subseteq S_0$. When this cofinal subfamily is countable, $\mathcal{S}$ is said to have a \textit{countable base}. A bornology $\mathcal{S}$ is said to have a \textit{compact (totally bounded) base} if it has a base consisting of compact (totally bounded) sets.    
 
 For $A \subseteq X$, the \textit{$\epsilon$-enlargement of $A$}, denoted by $B_d(A, \epsilon)$, is defined as, $B_d(A, \epsilon) = \{x \in X: d(x, A) < \epsilon\}$.

The Attouch-Wets convergence and the Hausdorff distance convergence can also be deduced as particular cases of another set convergence, called \textit{bornological convergence}, introduced by Lechicki, Levi, and Spakowski (\cite{BCL}).  The bornological convergence corresponding to a bornology $\mathcal{S}$ is denoted by $\mathcal{S}$-convergence. In spite of the fact that $\mathcal{S}$-convergence need not be topological in general, it has gained the interest of several researchers (\cite{beer2009operator,beer2023bornologies,Bcas,gapexcess,PbciAWc,Ucucas,Idealtopo}). When $\mathcal{S}$-convergence is topological on $CL(X)$, we denote the topology of $\mathcal{S}$-convergence by $\tau(\mathcal{S})$.  
 
 
 \begin{definition}\normalfont(\cite{BCL})
 	Suppose $\mathcal{S}$ is an ideal of subsets of a metric space $(X,d)$. A net $(A_\lambda)$ is said to be $\mathcal{S}$-convergent to $A$ in $\mathcal{P}_0(X)$ if for each $\epsilon > 0$ and $S \in \mathcal{S}$, eventually, we have $A_\lambda \cap S \subseteq B_d(A, \epsilon)\ \ \text{and}\ \ A \cap S \subseteq B_d(A_\lambda,\epsilon)$.
 \end{definition}
 
 In special cases, when $\mathcal{S} = \mathcal{B}_d(X)$ and $\mathcal{P}_0(X)$ the $\mathcal{S}$-convergence reduces to the Attouch-Wets convergence and Hausdorff distance convergence, respectively.



\begin{definition}\normalfont(\cite{Idealtopo,cao})
	Let $(X,d)$ be a metric space and $\mathcal{S} \subseteq \mathcal{P}_0(X)$. The topology $\tau_{\mathcal{S},d}$ on $\mathcal{P}_0(X)$ is determined by the uniformity $\mathcal{U}_{\mathcal{S},d}$ having base consisting of sets of the form $$ [S,\epsilon] = \{(A,C) \in \mathcal{P}_0(X) \times \mathcal{P}_0(X) : \sup_{x \in S}\vert d(x,A) - d(x,C)\vert < \epsilon\},$$ where $S \in \mathcal{S}$ and $\epsilon >0$.
	Clearly, a net $(A_\lambda)$ is $\tau_{\mathcal{S},d}$-convergent to $A$ if and only if  for any $S \in \mathcal{S}$ and $\epsilon > 0$, eventually, $\vert d(x,A) - d(x, A_\lambda)\vert < \epsilon$ for all $x \in S$.
	
%
\end{definition}	


For $\mathcal{S} = \mathcal{F}(X)$, $\mathcal{B}_d(X)$, and $\mathcal{P}_0(X)$, the topology $\tau_{\mathcal{S},d}$ reduces to the Wijsman, Attouch-Wets, and the classical Hausdorff distance topologies, respectively. See (\cite{francaviglia1985quasi,wijsconvergence}) for more on the Wijsman topology. We will use the following folklore in the sequel.

%
 
 \begin{proposition}\label{wijsmancoincidence}
 	Let $(X,d)$ be a metric space and let $\mathcal{S}$ be a bornology on $X$. Then $\tau_{\mathcal{S},d} = \tau_{W_{d}}$ on $CL(X)$ if and only if $\mathcal{S} \subseteq \mathcal{T}\mathcal{B}_{d}(X)$.
 \end{proposition} 
 \begin{proof}
 	Suppose $\tau_{\mathcal{S},d} = \tau_{W_{d}}$. Let $\mathcal{S} \nsubseteq  \mathcal{T}\mathcal{B}_d(X).$ Then there exists an $S \in \mathcal{S}$ such that $S$ is not totally bounded. Choose $r > 0$ such that for every $F \in \mathcal{F}(X)$, there exists $x_{F} \in S\setminus B_d(F,r)$. Define $ A_{F} = F$ for each $F \in \mathcal{F}(X)$ and $A = X$. Direct $\mathcal{F}(X)$ by set inclusion. Then $(A_{F})$ is a net in $CL(X)$. We claim that the net $(A_{F})$ is $\tau_{W_{d}}$-convergent to $A$ while $\tau_{\mathcal{S},d}$-convergence fails. Take $F_0 \in \mathcal{F}(X)$ and $\epsilon > 0$. Then for any $F \supseteq F_0$, $d(x,A_{F}) = 0 \quad \forall x \in F_0$. So for any $F \supseteq F_0$, we have $\displaystyle{\sup_{x \in F}}\vert d(x,A_F) - d(x,A)\vert = 0$. Thus, the net $(A_{F})$ is $\tau_{W_{d}}$-convergent to $A$. However, for any $F \in \mathcal{F}(X)$ there is an $x_{F} \in S\setminus B_d(F,r)$, which gives $d(x_{F},A_{F}) \geq r,$ and thus, $\displaystyle{\sup_{x \in S}} \vert d(x,A_{F}) - d(x,A) \vert > r$. Hence, the $\tau_{\mathcal{S},d}$-convergence of net $(A_{F})$ to $A$ fails. This contradicts the hypothesis.\\	
 	 Conversely, let $\mathcal{S} \subseteq  \mathcal{T}\mathcal{B}_{d}(X)$. We show that $\tau_{\mathcal{S},d} \subseteq \tau_{W_{d}}.$ Take any $S \in \mathcal{S}$ and $\epsilon > 0$. Then by the hypothesis, there exists $F \in \mathcal{F}(X)$ such that $S \subseteq B_d(F,{\frac{\epsilon}{4}})$. Let $(A,C) \in [F, \frac{\epsilon}{4}]$. Then $\displaystyle{\sup_{x \in F}}\vert d(x,A) - d(x,C)\vert < \frac{\epsilon}{4}$. Pick any $x \in S$. Choose $y \in F$ such that $d(x,y) < \frac{\epsilon}{4}$. Then \begin{eqnarray*}
 	 d(x,A) \leq d(x,y) + d(y,A) < d(y,A) + \frac{\epsilon}{4}\\ \text{and}~ d(y,C) \leq d(x,y) + d(x,C) < d(x,C) + \frac{\epsilon}{4}.\end{eqnarray*} So we get, \begin{eqnarray*}
 		d(x,A) - d(x,C) \leq d(y,A) - d(y,C) + \frac{2 \epsilon}{4} < \frac{3 \epsilon}{4}.\end{eqnarray*} Similarly, $d(x,C) - d(x,A) < \frac{3 \epsilon}{4}$. Thus, $\displaystyle{\sup_{x \in S}} \vert d(x,A) - d(x,C)\vert < \epsilon$. So $(A,C) \in [S, \epsilon]$.  
 \end{proof}
 
	\section{ A New class of Bornologies}

The properties of a set convergence which is determined by a bornology $\mathcal{S}$ may be characterized in terms of the properties satisfied by the  bornology $\mathcal{S}$. For example, the bornological convergence is topological on $CL(X)$ whenever the  underlying bornology is stable under small enlargements (\cite{BCL}). Similarly, if the bornology $\mathcal{S}$ has a countable base, then the topology $\tau_{\mathcal{S},d}$ is metrizable (\cite{Idealtopo}). To study various topological properties of the space $(CL(X),\tau_{\mathcal{S},d})$, we define and study in this section two properties for an arbitrary family of subsets of a metric space.

\begin{definition}\normalfont\label{Sseperable}
	Let $(X,d)$ be a metric space and let $\mathcal{S} \subseteq \mathcal{P}_0(X)$. Then we say \textit{$X$ is $\mathcal{S}$-separable} if there is a countable subcollection $\{S_n: n \in \mathbb{N}\}$ of $\mathcal{S}$ such that $\displaystyle{X = \overline{\bigcup_{n \in \mathbb{N}}S_n}}$. 
\end{definition}
If $(X,d)$ is separable and $\mathcal{S} \subseteq \mathcal{P}_0(X)$ forms a cover of $X$, then $X$ is $\mathcal{S}$-separable. However, the converse fails as $X$ is always $\mathcal{P}_0(X)$-separable but $(X,d)$ need not be separable.


Recall that a topological space $(X, \tau)$ is called \textit{$\sigma$-compact} if it is a countable union of compact subsets. If $X$ contains
a dense $\sigma$-compact subspace, then it is called \textit{almost $\sigma$-compact}. Clearly, for a metric space $(X,d)$, $X$ is $\mathcal{K}(X)$-separable if and only if $X$ is almost $\sigma$-compact, and $X$ is $\mathcal{F}(X)$-separable if and only if $X$ is separable.     

\begin{proposition}\label{TBsepimpliessep}
Let $(X,d)$ be a metric space and let $\mathcal{S}$ be a bornology such that $\mathcal{S} \subseteq \mathcal{T}\mathcal{B}_d(X)$. Then $X$ is $\mathcal{S}$-separable if and only if $X$ is separable.	 
\end{proposition}
\begin{proof}
It is enough to show that if $X$ is $\mathcal{S}$-separable, then $X$ is separable. Let $\{S_n: n \in \mathbb{N}\}\subseteq \mathcal{S}$ be such that $\displaystyle{X = \overline{\bigcup_{n \in \mathbb{N}}S_n}}$. For each $m, n \in \mathbb{N}$, choose $F_{m,n} \in \mathcal{F}(X)$ such that $F_{m,n} \subseteq S_n \subseteq B_d(F_{m,n}, \frac{1}{m})$. Set $A_n = \cup\{F_{m,n}: m \in \mathbb{N}\}$. Clearly, for each $n \in \mathbb{N}$, $A_n$ is countable. Put $\displaystyle{A = \bigcup_{n \in \mathbb{N}}A_n}$. We show that $A$ is dense in $(X,d)$. Suppose $x \in X$ and $\epsilon > 0$. Choose $m_0 \in \mathbb{N}$ such that $\frac{2}{m_0} < \epsilon$. Then for some $n_0 \in \mathbb{N}$, there exists $y \in S_{n_0}$ such that $d(x,y) < \frac{1}{m_0}$. Consequently, there is a $z \in F_{m_0,n_0}$ satisfying $d(y,z) < \frac{1}{m_0}$. So $d(x,z) <\epsilon$.       
\end{proof}
	
%

\begin{definition}\normalfont\label{totally bounded family}
	Let $(X,d)$ be a metric space. A nonempty family $\mathcal{S} \subseteq \mathcal{P}_0(X)$ is said to be a \textit{totally bounded family on $X$} (or \textit{totally bounded family on $(X,d)$}), if there is a countable subcollection $\{S_n: n \in \mathbb{N}\}$ of $\mathcal{S}$ such that for each $S \in \mathcal{S}$ and $\epsilon > 0$, we have $\displaystyle{S \subseteq \bigcup_{n \in F}B_d(S_n, \epsilon)}$ for some finite subset $F$ of $\mathbb{N}$.
	  
	Note that whenever $\mathcal{S}$ is a bornology the finite set $F$ may be taken as singleton set as one can replace the collection $\{S_n: n \in \mathbb{N}\}$ with $\displaystyle{\{\bigcup_{n \in F}S_n : F \in \mathcal{F}(\mathbb{N})\}}$. 
\end{definition}

It is easy to see that every bornology $\mathcal{S}$ with a countable base forms a totally bounded family on $(X,d)$. However, the converse need not be true as shown by the following example.

\begin{example}\label{Not countable base}
	\normalfont Let $(\mathbb{R},d_u)$ represent the set of real numbers with the usual metric. Then the bornology $\mathcal{F}(\mathbb{R})$ does not have a countable base. However, $\mathcal{S}(\mathbb{Q}) = \{\{a\} \vert \  \ a \in \mathbb{Q}\}$ is a countable subset of $\mathcal{F}(\mathbb{R})$ such that $\mathcal{F}(\mathbb{R})$ is a totally bounded family on $\mathbb{R}$. To see this, let $F = \{x_1, \ldots, x_n\} \in \mathcal{F}(\mathbb{R})$ and $\epsilon > 0$. Then for each $x_i \in F$ there is an $a_i \in \mathbb{Q}$ such that $d_u(x_i, a_i) <  \epsilon$. Consequently, $\displaystyle{F \subseteq \bigcup_{i=1}^nB_{d_u}(\{a_i\}, \epsilon)}$.\qed 	 
\end{example}
   

\begin{proposition}\label{P2impliesP1}
	Let $(X,d)$ be a metric space and let $\mathcal{S}$ be a bornology on $X$. If $\mathcal{S}$ is a totally bounded family on $X$, then $X$ is $\mathcal{S}$-separable.   
\end{proposition}

\begin{proof}
	Let $\{S_n: n \in \mathbb{N}\} \subseteq \mathcal{S}$ be such that $S_n \subseteq S_{n+1}$ and for any $S \in \mathcal{S}$ and $\epsilon > 0$ there is an $N \in \mathbb{N}$ with $S \subseteq B_d(S_N, \epsilon)$. We claim that $\displaystyle{X = \overline{\bigcup_{n \in \mathbb{N}}S_n}}$. Pick $x \in X$ and $\epsilon > 0$. Since $\{x\} \in \mathcal{S}$, there is an $n_0 \in \mathbb{N}$ with $\{x\} \subseteq B_d(S_{n_0}, \epsilon)$. So $\displaystyle{B_d(x, \epsilon) \cap (\bigcup_{n \in \mathbb{N}}S_n) \neq \emptyset}$. Thus, $X$ is $\mathcal{S}$-separable.  
\end{proof}

The converse of Proposition \ref{P2impliesP1} need not hold in general. 

\begin{example}\label{counterexample1}
	Suppose $X = \mathbb{R}$ and $d = d_{0,1}$, the $0$ - $1$ discrete metric. Consider the bornology $\mathcal{B}(\mathbb{R})$ having a base, $\{[-n,n] \cup C: C \in \mathcal{C}(\mathbb{Q}^c),~ n \in \mathbb{N}\}$, where $\mathcal{C}(\mathbb{Q}^c)$ denotes the collection of all nonempty countable subsets of $\mathbb{Q}^c$. 
	Then $X$ is $\mathcal{B}(\mathbb{R})$-separable, since $\displaystyle{\mathbb{R} = \bigcup_{n \in \mathbb{N}}[-n,n]}$. We show that $\mathcal{B}(\mathbb{R})$ is not a totally bounded family on $\mathbb{R}$. For this, we prove that for any countable subcollection $\{B_n: n \in \mathbb{N}\}$ of $\mathcal{B}(\mathbb{R})$ there exist $B \in \mathcal{B}(\mathbb{R})$ and $\epsilon > 0$ such that $B \nsubseteq B_d(B_n, \epsilon)$ for any $n \in \mathbb{N}$. It is sufficient to prove the claim when each $B_n$ is of the form $[-n,n]\cup C_n$, where $C_n \in \mathcal{C}(\mathbb{Q}^c)$. Let $\displaystyle{D = \mathbb{Q}^c \setminus \left(\bigcup_{n \in \mathbb{N}}C_n\right)}$. Clearly, $D$ is an uncountable subset of $\mathbb{Q}^c$. Furthermore, for each $n \in \mathbb{N}$, $(n, \infty) \cap D \neq \emptyset$, and $(-\infty, -n) \cap D \neq \emptyset$. Consequently, we can find $a_n, b_n \in D$ such that $a_n < -n, b_n > n$. Let $B = \{a_n: n \in \mathbb{N}\} \cup \{b_n: n \in \mathbb{N}\}$. Then $B \in \mathcal{B}(\mathbb{R})$ and for $0 < \epsilon < 1$, $B \nsubseteq B_{d_{0,1}}(B_n, \epsilon)$ for any $n \in \mathbb{N}$. Hence $\mathcal{B}(\mathbb{R})$ is not a totally bounded family on $\mathbb{R}$.\qed    
\end{example} 

\begin{example}\label{counterexample2}
	Suppose $(X,d)= (\mathbb{R}^2, d_s)$, where $$d_s(x,y) = \begin{cases}
		||x-y||_2 & \text{if $x = \alpha y$} \\
		&       \text{for some $\alpha \in \mathbb{R};$}\\
		||x||_2 + ||y||_2 & \text{otherwise}
	\end{cases}$$
	and $\parallel \cdot \parallel_2$ denotes the Euclidean norm on $\mathbb{R}^2$. Then $d_s$ is a metric known as the \textit{French Metro metric} on $\mathbb{R}^2$ (\cite{deza2009encyclopedia}).
	For each $n \in \mathbb{N}$, define $S_n = \{x \in \mathbb{R}^2: d_s(x,0)= n\}$ and $A_n = \{x \in \mathbb{R}^2 \setminus (\mathbb{Q} \times \mathbb{Q}): d_s(x,0)\leq n \}$. Clearly, $\displaystyle{\bigcup_{n \in \mathbb{N}}A_n = \mathbb{R}^2\setminus (\mathbb{Q} \times \mathbb{Q})}$. Let $\displaystyle{S = \bigcup_{n \in \mathbb{N}}S_n}$ and $\mathcal{C}(S)$ be the collection of all nonempty countable subsets of $S$. Set $\mathcal{B}_0 = \{A_n \cup A \cup F : A \in \mathcal{C}(S), F \in \mathcal{F}(\mathbb{R}^2), n \in \mathbb{N}\}$. Then $\mathcal{B}_0$ forms a base for a bornology on $\mathbb{R}^2$, say $\mathcal{B}(\mathbb{R}^2)$.
	 Take $\displaystyle{(a,b) \in \mathbb{R}^2\setminus \left(\bigcup_{n \in \mathbb{N}}A_n\right)}$ and $\epsilon > 0$. If $(a,b) = (0,0)$, then $B_{d_s}((0,0), \epsilon) = \{x \in \mathbb{R}^2:~ \parallel x \parallel_2 < \epsilon\}$. Clearly, $B_{d_s}((0,0), \epsilon) \cap (\displaystyle{\bigcup_{n \in \mathbb{N}}A_n}) \neq \emptyset$. So $\displaystyle{(0,0) \in \overline{\bigcup_{n \in \mathbb{N}}A_n}}$. Otherwise when $(a,b) \neq (0,0)$, choose $\alpha \in (1 - \frac{\epsilon}{\sqrt{a^2 + b^2}}, 1 + \frac{\epsilon}{\sqrt{a^2 + b^2}} )$ such that $(\alpha a, \alpha b ) \in \mathbb{R}^2\setminus (\mathbb{Q} \times \mathbb{Q})$. Consequently, $(\alpha a , \alpha b)\in B_{d_s}((a,b), \epsilon)$. Therefore $\displaystyle{\mathbb{R}^2 = \overline{\bigcup_{n \in \mathbb{N}}A_n}}$. Thus, $(\mathbb{R}^2, d_s)$ is $\mathcal{B}(\mathbb{R}^2)$-separable.
	
	We next show that for any countable subcollection $\{B_n: n \in \mathbb{N}\}$ of $\mathcal{B}_0$, there exist $B \in \mathcal{B}(\mathbb{R}^2)$ and $\epsilon > 0$ such that for any $n \in \mathbb{N}$, we have $B \nsubseteq B_d(B_n, \epsilon)$. Thus, $\mathcal{B}(\mathbb{R}^2)$ is not a totally bounded family on $\mathbb{R}^2$. 
	
	To see this, let each $B_n$ be of the form $A_n \cup C_n \cup F_n$, where $F_n \in \mathcal{F}(\mathbb{R}^2)$ and $C_n \in \mathcal{C}(S)$. Define $\displaystyle{C = \mathbb{R}^2\setminus \left(\bigcup_{n \in \mathbb{N}}C_n\right)}$. Since $\displaystyle{\bigcup_{n \in \mathbb{N}}C_n}$ is a countable subset of $S$ and each $S_n$ is uncountable, there exists an $x_n \in S_n \cap C$ for each $n \in \mathbb{N}$.  Let $B = \{x_n: n\ \in \mathbb{N}\}$ and $\epsilon = \frac{1}{2}$.  Note that $B_{d_{s}}(A_n, \epsilon)$ and $B_{d_{s}}(F_n, \epsilon)$ are bounded, then by construction of $B$, we can find an $x_m \in B$ such that $x_m \notin B_{d_{s}}(A_n, \epsilon) \cup B_{d_{s}}(F_n, \epsilon) \cup C_n$. Observe that for each $n \in \mathbb{N}$, $\{B_{d_{s}}(x,n): x \in S_n \}$ form a pairwise disjoint collection of open balls. Also if $x \in S_n, y \in S_m$ such that $n \neq m$, we have $d_s(x,y) \geq 1$.  Consequently, $x_m  \notin B_{d_{s}}(C_n, \epsilon)$. Thus, $B \nsubseteq B_{d_{s}}(B_n, \epsilon)$.\qed   
\end{example}

We now give a sufficient condition for when the converse of Proposition \ref{P2impliesP1} is true. 

\begin{proposition}\label{Equivalence condition}
	Let $(X,d)$ be a metric space and let $\mathcal{S}$ be a bornology on $X$ with a compact base. Then \textit{$X$ is $\mathcal{S}$-separable} if and only if $\mathcal{S}$ is a totally bounded family on $X$.
\end{proposition}

\begin{proof}
	We only need to show that if $X$ is $\mathcal{S}$-separable, then $\mathcal{S}$ is a totally bounded family on $X$. Take $S \in \mathcal{S}$ and $\epsilon > 0$. Since $\mathcal{S}$ has a compact base, there exists a compact set $K \in \mathcal{S}$ such that $S \subseteq K$. Let $\{S_n: n \in \mathbb{N}\} \subseteq \mathcal{S}$ such that $\displaystyle{X = \overline{\bigcup_{n \in \mathbb{N}}S_n}}$. For each $x \in K$,  there exists $n_x \in \mathbb{N}$ such that $x \in B_d(S_{n_{x}}, \epsilon)$, which gives, $\displaystyle{K \subseteq \bigcup_{x \in K}B_d(S_{n_{x}}, \epsilon)}$. So there are $ x_1,x_2,\ldots,x_{k} \in K$ such that $\displaystyle{K \subseteq \bigcup_{i=1}^{k}B_d(S_{n_{x_{i}}}, \epsilon)}$. Consequently, $S \subseteq \displaystyle{\bigcup_{i=1}^{k} B_d(S_{n_{x_{i}}}, \epsilon)}$. \end{proof}

The following example shows that the condition of $\mathcal{S}$ having a compact base need not be necessary in Proposition \ref{Equivalence condition}. In fact, one can easily show that Proposition \ref{Equivalence condition} holds for a bornology $\mathcal{S}$ with a totally bounded base. 

\begin{example}
	Suppose $X = \{x = (x_n): x_n = 0 \text{ eventually}\}$ with metric $d$ defined as, $\displaystyle{d((x_n),(y_n)) = \sup_{n \in \mathbb{N}}\vert x_n - y_n \vert}$. Consider the bornology $\mathcal{B}_d(X)$ on $X$. Since $\mathcal{B}_d(X)$ has a countable base,  $\mathcal{B}_d(X)$ is a totally bounded family on $X$. However, $\mathcal{B}_d(X)$ does not have a compact base, as $X$ is infinite dimensional.\qed
\end{example}


\begin{definition}\normalfont (\cite{beer2009totalboundedness, BCL}) Let $(X,d)$ be a metric space  and let $\mathcal{S} \subseteq \mathcal{P}_0(X)$.
A nonempty subset $A$ of $X$ is called \textit{weakly $\mathcal{S}$-totally bounded} if for each $\epsilon > 0$ there exist $S_1,\ldots,S_n \in \mathcal{S}$ such that $\displaystyle{A \subseteq  \bigcup_{i=1}^{n}B_d(S_i, \epsilon)}$. Note that whenever $\mathcal{S}$ is a bornology on $X$, then for each $\epsilon > 0$ we can find a $S \in \mathcal{S}$ such that $A \subseteq B_d(S, \epsilon)$. 
\end{definition} 
Denote by $\mathcal{S}^*$ the collection of all weakly $\mathcal{S}$-totally bounded nonempty subsets of $X$. Clearly, $\mathcal{S} \subseteq \mathcal{S}^*$. Also, it is easy to verify that $\mathcal{S}^*$ forms a bornology on $X$. 
The next result shows that the condition of bornology $\mathcal{S}$ being a totally bounded family on $X$ is weaker than that of $X$ is weakly $\mathcal{S}$-totally bounded.
\begin{proposition}
	Let $(X,d)$ be a metric space and let $\mathcal{S}$ be a bornology on $X$. If $X$ is weakly $\mathcal{S}$-totally bounded, then $\mathcal{S}$ is a totally bounded family on $X$. 
\end{proposition}
\begin{proof}
	Let $S \in \mathcal{S}$ and $\epsilon > 0$. For each $n \in \mathbb{N}$, choose a finite subset $F_n$ of $\mathcal{S}$ such that $\displaystyle{X = \bigcup_{S'\in F_n}B_d(S',1/ n)}$. Define $\beta =\cup\{F_n: n \in \mathbb{N} \}$. Then $\beta$ is a countable subfamily of $\mathcal{S}$. Choose $n_0 \in \mathbb{N}$ such that $\frac{1}{n_{0}} < \epsilon$ and $\displaystyle{S \subseteq \bigcup_{S'\in F_{n_0}}B_d(S', 1/n_{0}) \subseteq \bigcup_{S'\in F_{n_0}}B_d(S', \epsilon)}$.    
\end{proof}
\begin{proposition}\label{weaklytotallybonded}
	Let $(X,d)$ be a metric space and let $\mathcal{S}$ be a bornology on $X$. Then $\mathcal{S}$ is a totally bounded family on $X$ if and only if $\mathcal{S}^*$ is a totally bounded family on $X$.
\end{proposition}
\begin{proof}
	Let $\{S_n: n \in \mathbb{N}\}$ be the subcollection of $\mathcal{S}$ satisfying Definition \ref{totally bounded family} for $\mathcal{S}$ such that $S_n \subseteq S_{n+1}$. Let $A \in \mathcal{S}^*$ and  $\epsilon > 0$. Then there exists $S \in \mathcal{S}$ such that $ A \subseteq B_d(S, \frac{\epsilon}{2})$. Consequently, there exists $n_0 \in \mathbb{N}$ such that $S \subseteq B_d(S_{n_0}, \frac{\epsilon}{2})$. Therefore, $A \subseteq B_d(S_{n_0}, \epsilon)$. Thus, $\mathcal{S}^*$ is a totally bounded family on $X$. 
	
	Conversely, suppose $\mathcal{S}^*$ is a totally bounded family on $X$. Choose $\{S_n: n \in \mathbb{N}\} \subseteq \mathcal{S}^*$ satisfying Definition \ref{totally bounded family} for $\mathcal{S}^*$ such that $S_n \subseteq S_{n+1}$. Since each $S_n$ is weakly $\mathcal{S}$-totally bounded, for each $n \in \mathbb{N}$ and $m \in \mathbb{N}$ we can  find $B_{n,m} \in \mathcal{S}$ such that $S_n \subseteq B_d(B_{n,m}, \frac{1}{m})$. Let $\mathscr{A}_{n,m} = \{B_{n,m}: m \in \mathbb{N}\}$ for $n \in \mathbb{N}$ and $\displaystyle{\mathscr{A} = \bigcup_{n \in \mathbb{N}}\mathscr{A}_{n,m}}$. Clearly, $\mathscr{A}$ is a countable subcollection of $\mathcal{S}$. Let $S \in \mathcal{S}$ and $\epsilon > 0$. Choose $m \in \mathbb{N}$ such that $\frac{2}{m} < \epsilon$. Since $S \in \mathcal{S}^*$, there exists an $n_0 \in \mathbb{N}$ such that $S \subseteq B_d(S_{n_0}, \frac{1}{m})$. Also there exists $B_{n_0, m} \in \mathcal{A}$ such that $S_{n_0} \subseteq B_d(B_{n_0,m}, \frac{1}{m})$. So $S \subseteq B_d(B_{n_0, m}, \epsilon)$. Thus, $\mathcal{S}$ is a totally bounded family on $X$.   
	\end{proof}
	
\section{A metric topology on $CL(X)$ and Submetrizabililty of $\tau_{\mathcal{S},d}$}
Recall that a topological space $(X,\tau)$ is called \textit{submetrizable} if $X$ admits a weaker metrizable topology. In this section, we make use of the properties studied in the previous section to characterize submetrizability of the space $(CL(X),\tau_{\mathcal{S}, d})$. First, we define a metric on the set $CL(X)$ for an $\mathcal{S}$-separable metric space $(X,d)$.     
\begin{proposition}\label{metric}
	Let $(X,d)$ be a metric space and let $\mathcal{S}$ be a bornology on $X$. Suppose there exists a countable subcollection $\mathcal{A} = \{S_i: i \in \mathbb{N}\}$ of $\mathcal{S}$ that satisfies Definition \ref{Sseperable} for $\mathcal{S}$, that is, $\displaystyle{X = \overline{\bigcup_{S_i \in \mathcal{A}}S_i}}$. 
	 Then the map $d_{\mathcal{S}}^{\mathcal{A}}:CL(X) \times CL(X) \rightarrow \mathbb{R}$ defined by $$d_{\mathcal{S}}^{\mathcal{A}}(A,C) = \sum_{i=1}^{\infty}2^{-i}\min\{1, \sup_{x \in S_{i}}\vert d(x,A) - d(x,C) \vert\},$$
	is a metric on $CL(X)$.
	\end{proposition}

\begin{proof}
	Let $A,C,$ and $D \in CL(X)$. Clearly, $d_{\mathcal{S}}^{\mathcal{A}}(A,A) = 0$. Suppose $d_{\mathcal{S}}^{\mathcal{A}}(A,C) = 0$. Then for each $i \in \mathbb{N}$, $d(x, A) = d(x,C)\quad \forall x \in S_i$. Since $d(\cdot, A)$ and $d(\cdot, C)$ are continuous functions on $X$ and they agree on the dense subset $\displaystyle{\bigcup_{i \in \mathbb{N}}S_i}$ of $X$,  $d(\cdot, A) = d(\cdot, C)$ on $X$. Consequently, $A = C$. The symmetric property is obvious from the definition of $d_{\mathcal{S}}^{\mathcal{A}}$. To see the triangle inequality, consider $x \in S_i$ for any $i \in \mathbb{N}$, then we have $\vert d(x,A) - d(x,C) \vert \leq \vert d(x,A) - d(x, D) \vert + \vert d(x,D) - d(x,C) \vert$ for $A, C, D \in CL(X)$. So $$\sup_{x \in S_{i}}\vert d(x,A) - d(x, C)\vert \leq \sup_{x \in S_{i}} \vert d(x,A) - d(x,D) \vert + \sup_{x \in S_{i}} \vert d(x,D) -d(x,C) \vert.$$ 
	Thus, $d_{\mathcal{S}}^{\mathcal{A}}(A,C) \leq d_{\mathcal{S}}^{\mathcal{A}}(A,D) + d_{\mathcal{S}}^{\mathcal{A}}(D,C)$.      
\end{proof}

\begin{corollary}
	Let $(X,d)$ be a metric space and let $\mathcal{S}$ be a bornology on $X.$ If $\mathcal{S}$ has a countable base $ \mathcal{A} = \{S_i : \  \ n \in \mathbb{N}\}$, then $$d_{\mathcal{S}}^{\mathcal{A}}(A,C) = \sum_{i=1}^{\infty} 2^{-i} \min \{1, \sup_{x \in S_{i}}|d(x, A) - d(x, C)|\}$$ is a metric on $CL(X)$. 
\end{corollary} 
The next proposition shows that the topology on $CL(X)$ corresponding to metric $d_{\mathcal{S}}^{\mathcal{A}}$ is always finer than the Wijsman topology, $\tau_{W_{d}}$ on $CL(X)$.  
\begin{proposition}\label{dSAfinerwijsman}
	Let $(X,d)$ be a metric space and let $\mathcal{S}$ be a bornology on $X$. Suppose $\mathcal{A} = \{S_i: i \in \mathbb{N}\}$ is a countable subcollection of $\mathcal{S}$ such that $\displaystyle{X = \overline{\bigcup_{S_i \in \mathcal{A}}S_i}}$, then the topology generated by metric $d_{\mathcal{S}}^{\mathcal{A}}$ is always finer than the Wijsman topology, $\tau_{W_{d}}$ on $CL(X)$. 
\end{proposition}
\begin{proof}
	Let $(A_\lambda)$ be a net that converges to $A$ in $(CL(X), d_{\mathcal{S}}^{\mathcal{A}})$. Take $0 < \epsilon <1$ and $x \in X$. Choose $i_0 \in \mathbb{N}$ such that $d(x,y) < \frac{\epsilon}{4}$ for some $y \in S_{i_0}$. By the assumption, there exists a $\lambda_0$ such that $d_{\mathcal{S}}^{\mathcal{A}}(A, A_{\lambda}) < \frac{\epsilon}{2^{i_0 +2}}$ $\forall \lambda \geq \lambda_0$. Then $\displaystyle{\sup_{z \in S_{i_0}}\vert d(z, A) - d(z, A_{\lambda})\vert < \frac{\epsilon}{2}}$ for  $\lambda \geq \lambda_0$. So for $\lambda \geq \lambda_0$, we have $d(x,A) \leq d(x,y) + d(y, A) < d(x,y) + d(y, A_{\lambda}) + \frac{\epsilon}{2}< d(x, A_{\lambda}) + \epsilon$. Consequently, $\vert d(x,A) - d(x, A_{\lambda})\vert < \epsilon$ $\forall \lambda \geq \lambda_0$. Therefore, $\tau_{W_{d}}$-convergence of $(A_{\lambda})$ to $A$ is established.      
\end{proof}
 
 Recall that a hyperspace topology $\tau$ on $CL(X)$ is said to be \textit{admissible} if the mapping $\psi: (X,d) \to (CL(X), \tau)$ defined by $\psi(x) = \{x\}$ for $x \in X$ is an embedding. 
\begin{theorem}
	Let $(X,d)$ be a metric  space and let $\mathcal{S}$ be a bornology on $X$. Suppose $\mathcal{A}=\{S_i: i \in \mathbb{N}\}$ is a countable  subcollection of $\mathcal{S}$ such that $\displaystyle{X = \overline{\bigcup_{i \in \mathbb{N}}S_i}}$. Then the topology generated by metric $d_{\mathcal{S}}^{\mathcal{A}}$ on $CL(X)$ is admissible. 
\end{theorem}
\begin{proof}
	 Suppose $(z_n)$ is a sequence in $(X,d)$ that converges to $z$ in $(X,d)$. Let $0 < \epsilon < 1$. Choose $n_0 \in \mathbb{N}$ such that $d(z, z_n) < \frac{\epsilon}{2}$ $\forall n \geq n_0$. So for any $x \in X$, we have $\vert d(x,z) - d(x,z_n)\vert < \frac{\epsilon}{2}$ $\forall n \geq n_0$. Consequently, $d_{\mathcal{S}}^{\mathcal{A}}(\{z\}, \{z_n\}) < \epsilon$ for $n \geq n_0$. Therefore, the continuity of $\psi$ follows. \\
	Next, we show the continuity of $\psi^{-1}:(\psi(X), d_{\mathcal{S}}^{\mathcal{A}}) \to (X,d)$. To see this, let $(\{y_n\})$ be a sequence that converges to $\{y\}$ with respect to $d_{\mathcal{S}}^{\mathcal{A}}$. Take $0 < \epsilon < 1$. Choose $i_0 \in \mathbb{N}$ such that $y \in B_d(S_{i_0}, \frac{\epsilon}{8})$. By the assumption, we have $n_0 \in \mathbb{N}$ such that $$\sum_{i =1}^{\infty}2^{-i}\min\{1, \sup_{x \in S_i}\vert d(x,\{y\}) - d(x,\{y_n\})\vert\}< \frac{\epsilon}{2^{i_0+4}} \quad \forall n \geq n_0.$$
	Consequently, $\displaystyle{\sup_{x\in S_{i_0}}\vert d(x,y) - d(x,y_n)\vert < \frac{\epsilon}{8}}$ $\forall n \geq n_0$. Choose $x_0 \in S_{i_0}$ with $d(x_0,y)< \frac{\epsilon}{8}$. Then $d(y, y_n) \leq d(y,x_0) + d(x_0, y_n)< 2d(x_0,y) + \frac{\epsilon}{8} < \epsilon$ $\forall n \geq n_0$. Therefore, for $n \geq n_0$, we have $d(y, y_n) < \epsilon$. Hence, $(y_n)$ converges to $y$ in $(X,d)$.
\end{proof}

The next example shows that the metric $d_{\mathcal{S}}^\mathcal{A}$ in Proposition \ref{metric} depends upon the choice of countable subcollection $\mathcal{A}$. In fact, corresponding to two different subcollections $\mathcal{A}_1$ and $\mathcal{A}_2$ from $\mathcal{S}$ satisfying Definition \ref{Sseperable}, the resultant metrics $d_{\mathcal{S}}^{\mathcal{A}_1}$ and $d_{\mathcal{S}}^{\mathcal{A}_2}$ on $CL(X)$ may not be necessarily equivalent. Recall that if $(X,d)$ is separable and  $\mathcal{S} = \mathcal{F}(X)$, then for any countable dense set $\{x_n: n \in \mathbb{N}\}$ of $X$, the metric $d_{\mathcal{F}(X)}^{\mathcal{A}}$, where $\mathcal{A} = \{\{x_n\}: n \in \mathbb{N}\}$ is compatible with the Wijsman topology $\tau_{W_d}$ on $CL(X)$ (see, Exercise 2.1.6 in \cite{ToCCoS}). We use this fact in our next example.

\begin{example}\label{metric_diffcollection}
	Let $(X,d) = (\mathbb{R}^2,d_e)$ be the Euclidean plane and let $\mathcal{S} = \{S \subseteq (\mathbb{N} \times \{0\}) \cup F: F \subseteq \mathbb{R}^2 \text{ is finite}\}$. It can be verified that the collections $\mathcal{A}_1 = \{\{(x,y)\}: x,y \in \mathbb{Q}\}$ and $\mathcal{A}_2 = \{\mathbb{N} \times \{0\}\} \cup \{\{(x,y)\}: x,y \in \mathbb{Q}\}$ satisfy Definition \ref{Sseperable}. Therefore, by Proposition \ref{metric}, $d_{\mathcal{S}}^{\mathcal{A}_1}$ and $d_{\mathcal{S}}^{\mathcal{A}_2}$ are two metrics on $CL(X)$. We show that these two metrics are not equivalent on $CL(X)$. Let for each $n \in \mathbb{N}$, $A_n$ denote the line $y = \frac{x}{n}$ and let $A$ be the horizontal axis. Then we show that the sequence $(A_n)$ does not converge to $A$ with respect to the metric $d_{\mathcal{S}}^{\mathcal{A}_2}$ while it converges to $A$ with respect to $d_{\mathcal{S}}^{\mathcal{A}_1}$.
	
	First we show convergence with respect to $d_{\mathcal{S}}^{\mathcal{A}_1}$. By the above discussion, it is enough to show that the sequence $(A_n)$ converges to $A$ in $(CL(X), \tau_{AW_{d}})$. 
	To see this, let $B$ be any bounded subset of $\mathbb{R}^2$ and $\epsilon > 0$. Then for some $r > 0$, we have $B \subseteq B_d((0, 0), r)$. Choose $n_0 \in \mathbb{N}$ such that $\frac{r}{n_0} < \epsilon$.  Pick $(x,0) \in A \cap B_d((0, 0), r)$. Then $d((x,0), A_n) \leq \frac{\vert x \vert}{n} <  \frac{r}{n}$. So for $n \geq n_0$, we have $d((x,0), A_n) < \epsilon$. Therefore, $A \cap B_d((0, 0), r) \subseteq B_d(A_n, \epsilon)$ for $n \geq n_0$. Similarly, let $ (x, \frac{x}{n}) \in A_n \cap B_d((0,0), r)$ for any $n \in \mathbb{N}$. Then $d((x, \frac{x}{n}), A)\leq \frac{\vert x\vert}{n} <  \frac{r}{n}$. Consequently, $A_n \cap B_d((0, 0), r) \subseteq B_d(A, \epsilon)$ for $n \geq n_0$. Thus, $(A_n)$ converges to $A$ in $(CL(X), \tau_{AW_{d}})$. 
	
	Next, we show that $(A_n)$ fails to converge to $A$ with respect to $d_{\mathcal{S}}^{\mathcal{A}_2}$. Let $\{(x_m,y_m): m \in \mathbb{N}\}$ be an enumeration for $\mathbb{Q}\times \mathbb{Q}$. So $\mathcal{A}_2 = \{\mathbb{N} \times \{0\}\} \cup \{\{(x_m,y_m)\}: m \in \mathbb{N}\}$.  Then
	\begin{eqnarray*}
		d_{\mathcal{S}}^{\mathcal{A}_2}(A, A_n) = \frac{1}{2}\min \{1, \sup_{(x,0) \in \{\mathbb{N} \times \{0\}\}}\vert d((x,0), A) - d((x,0), A_n)\vert\}\\ + \sum_{m =1}^{\infty}2^{-{(m+1)}} \min \{1, \vert d((x_m,y_m), A) - d((x_m,y_m), A_n)\vert \}  \geq \frac{1}{2}.
	\end{eqnarray*} 
	Thus, $(A_n)$ does not converges to $A$ with respect to $d_{\mathcal{S}}^{\mathcal{A}_2}$.\qed
\end{example}
\begin{proposition}\label{dSAequivalentmetric}
	Let $(X,d)$ be a metric space and let $\mathcal{S}$ be a bornology on $X$. Suppose $\mathcal{A} =\{S_n: n \in \mathbb{N}\}$ is a countable subcollection of $\mathcal{S}$ such that $\displaystyle{X = \overline{\bigcup_{S_n \in \mathcal{A}}S_n}}$ and $\mathcal{A}' = \{S'_n: n\in \mathbb{N}\}$, where $\displaystyle{S'_n = \bigcup_{i =1}^nS_i}$, then the metrics $d_{\mathcal{S}}^{\mathcal{A}}$ and $d_{\mathcal{S}}^{\mathcal{A}'}$  are equivalent on $CL(X)$.
\end{proposition}
\begin{proof}
	It is easy to verify that the metric $d_{\mathcal{S}}^{\mathcal{A}'}$ is stronger than the metric $d_{\mathcal{S}}^{\mathcal{A}}$. We show that the converse also holds. Let $(A_m)$ be a sequence that converges to $A$ in $(CL(X), d_{\mathcal{S}}^{\mathcal{A}})$. Take $\epsilon > 0$. Choose $0<r<1$ and $n_0 \in \mathbb{N}$ such that $r\left(1 - \frac{1}{2^{n_0}}\right)+ \frac{1}{2^{n_0}} < \epsilon$. Then there exists $m_0 \in \mathbb{N}$ such that $d_{\mathcal{S}}^{\mathcal{A}}(A_m, A)< \frac{r}{2^{n_0+1}}$ $\forall m \geq m_0$. So for $m \geq m_0$, we have $\displaystyle{\sup_{x \in S_n}\vert d(x, A_m) - d(x, A)\vert} < \frac{r}{2}$ for $n \leq n_0$. Consequently, $\displaystyle{\sup_{x \in S'_n}\vert d(x, A_m) - d(x, A)\vert} < r$ for $n \leq n_0$.
Therefore, for all $m \geq m_0$, we have \begin{eqnarray*}
d_{\mathcal{S}}^{\mathcal{A}'}(A_m, A) = \sum_{n=1}^{n_0}2^{-n}\min\{1, \sup_{x \in S'_n}\vert d(x, A_m)-d(x, A)\vert\}+\\ \sum_{n_0+1}^{\infty}2^{-n}\min\{1, \sup_{x \in S'_n}\vert d(x, A_m)-d(x, A)\vert\}\\< r\left(1 - \frac{1}{2^{n_0}}\right)+ \frac{1}{2^{n_0}} < \epsilon.
\end{eqnarray*}\end{proof}
	In (\cite{Idealtopo}), the authors showed that for two bornologies $\mathcal{S}$ and $\mathcal{C}$ on $(X,d)$, the topologies $\tau_{\mathcal{S},d}$ and $\tau_{\mathcal{C},d}$ are equal on $\mathcal{P}_0(X)$ if and only if $\mathcal{S}^* = \mathcal{C}^*$. The following theorem explores when the metrics $d_{\mathcal{S}}^{\mathcal{A}_1}$ and $d_{\mathcal{C}}^{\mathcal{A}_2}$ are equivalent on $CL(X)$, where $\mathcal{A}_1$ and $\mathcal{A}_2$ are countable subcollections of $\mathcal{S}$ and $\mathcal{C}$, respectively such that $\displaystyle{X = \overline{\bigcup_{S \in \mathcal{A}_1}S} = \overline{\bigcup_{C \in \mathcal{A}_2}C}}$. It is  easy to verify that if $\mathcal{A}_1 \subseteq \mathcal{A}_2^*$, then $\mathcal{A}_1^* \subseteq \mathcal{A}_2^*$. 
\begin{theorem}\label{metricdiffbornologies}
	Let $(X,d)$ be a metric space and let $\mathcal{S}$, $\mathcal{C}$ be two bornologies on $X$. Suppose $\mathcal{A}_1$ and $\mathcal{A}_2$ are countable subcollections of $\mathcal{S}$ and $\mathcal{C}$, respectively such that $\displaystyle{X = \overline{\bigcup_{S \in \mathcal{A}_1}S} = \overline{\bigcup_{C \in \mathcal{A}_2}C}}$. Then, the metrics $d_{\mathcal{S}}^{\mathcal{A}_1}$ and $d_{\mathcal{C}}^{\mathcal{A}_2}$ are equivalent  on $CL(X)$ if and only if $\mathcal{A}_1^* = \mathcal{A}_2^*$.   
\end{theorem}

\begin{proof}
	Suppose $\mathcal{A}_1^* = \mathcal{A}_2^*$. It is enough to show that every sequence converging in $(CL(X), d_{\mathcal{S}}^{\mathcal{A}_1})$ also converges in $(CL(X), d_{\mathcal{C}}^{\mathcal{A}_2})$. Let $(A_n)$ be a sequence in $CL(X)$ that converges to $A \in CL(X)$ with respect to $d_{\mathcal{S}}^{\mathcal{A}_1}$. By Proposition \ref{dSAequivalentmetric}, we can assume $\mathcal{A}_1 = \{S_i: i \in \mathbb{N}\}$, where $S_i \subseteq S_{i+1}$, and $\mathcal{A}_2 = \{C_i: i \in \mathbb{N}\}$, where $C_i \subseteq C_{i+1}$. Take $0< \epsilon < 1$. Choose $i_0 \in \mathbb{N}$ such that $\displaystyle{\frac{1}{2^{i_0}}< \frac{\epsilon}{8}}$. Then there exists $i_1 \in \mathbb{N}$ such that $i_{1} > i_{0}$ and $C_{i_0} \subseteq B_d(S_{i_1}, \frac{\epsilon}{8})$. Also there exists $n_0 \in \mathbb{N}$ such that $d_{\mathcal{S}}^{\mathcal{A}_1}(A, A_n) < \frac{\epsilon}{2^{i_1+3}} \quad \forall n \geq n_0$. Consequently, $\displaystyle{\sup_{y \in S_{i_1}}\vert d(y, A) - d(y, A_n)\vert < \frac{\epsilon}{8}}$ $\forall n \geq n_0$. Fix $n \geq n_0$ and let $x \in C_i$ for $i \leq i_0$. Since $C_{i_0} \subseteq B_d(S_{i_1}, \frac{\epsilon}{8})$, choose $y \in S_{i_1}$ such that $d(x,y) < \frac{\epsilon}{8}$. Then $d(x, A) \leq d(x,y) + d(y, A) < d(x,y) + d(y, A_n)+ \frac{\epsilon}{8}< 2d(x,y) + d(x, A_n) + \frac{\epsilon}{8}$. Consequently, $d(x, A) - d(x, A_n) < \frac{3\epsilon}{8}$. Similarly, $d(x, A_n) - d(x, A) < \frac{3\epsilon}{8}$. Therefore, $\displaystyle{\sup_{x \in C_i}\vert d(x, A) - d(x, A_n)\vert < \frac{\epsilon}{2}}$ for $i \leq i_0$. So \begin{eqnarray*}
		\sum_{i =1}^{i_0}2^{-i}\min\{1, \sup_{x \in C_i}\vert d(x, A) - d(x, A_n)\vert\} < \frac{\epsilon}{2}\left(1 - \frac{1}{2^{i_0}}\right). \text{ Also}\\
		\sum_{i =i_0+1}^{\infty}2^{-i}\min\{1, \sup_{x \in C_i}\vert d(x, A) - d(x, A_n)\vert\} \leq \frac{1}{2^{i_0}}.
	\end{eqnarray*} Hence for each $n \geq n_0$, we have $$d_{\mathcal{C}}^{\mathcal{A}_2}(A, A_n) = \sum_{i =1}^{\infty}2^{-i}\min\{1, \sup_{x \in C_i}\vert d(x, A) - d(x, A_n)\vert\} < \frac{\epsilon}{2}\left(1 - \frac{1}{2^{i_0}}\right) + \frac{\epsilon}{8} < \epsilon .$$\\
	Conversely, suppose $\mathcal{A}_1^* \neq \mathcal{A}_2^*$. Assume $\mathcal{A}_1^* \nsubseteq \mathcal{A}_2^*$. By the above discussion, $\mathcal{A}_1 \nsubseteq \mathcal{A}_2^*$. Then there is an $i_0 \in \mathbb{N}$ such that $S_{i_0} \notin \mathcal{A}_2^*$. For each $n \in \mathbb{N}$, define $A_n = \overline{C_n}$ and $A = X$. We claim that $(A_n)$ converges to $A$ with respect to $d_{\mathcal{C}}^{\mathcal{A}_2}$ while convergence to $A$ with respect to $d_{\mathcal{S}}^{\mathcal{A}_1}$ fails. Take $\epsilon > 0$. Choose $n_0 \in \mathbb{N}$ such that $\frac{1}{2^{n_0}} < \epsilon$. Consider \begin{eqnarray*}
		d_{\mathcal{C}}^{\mathcal{A}_2}(A, A_n) = \sum_{i =1}^{\infty}2^{-i}\min \{1, \sup_{x \in C_i}\vert d(x, X) - d(x, \overline{C_n})\vert\} \\ = \sum_{i = n+1}^{\infty}2^{-i} \min\{1, \sup_{x \in C_i}d(x, \overline{C_n})\} \leq \frac{1}{2^n}.
	\end{eqnarray*} So for $n \geq n_0$, $d_{\mathcal{C}}^{\mathcal{A}_2}(A, A_n) < \epsilon$. Since $S_{i_0} \notin \mathcal{A}_2^*$, there exists $r > 0$ such that for each $i' \in \mathbb{N}$, we have $S_{i_0} \nsubseteq B_d(C_{i'}, r)$. Consequently for each $n \in \mathbb{N}$, we have  \begin{eqnarray*}
		d_{\mathcal{S}}^{\mathcal{A}_1}(A,A_n) = \sum_{i =1}^{\infty}2^{-i}\min \{1, \sup_{x \in S_i} \vert d(x, X) - d(x, \overline{C_n})\vert\}\\ 
		= \sum_{i=1}^{\infty}2^{-i}\min \{1, \sup_{x \in S_i}d(x, \overline{C_n})\} \geq 2^{-i_0}\min\{1, \sup_{x \in S_{i_0}}d(x, \overline{C_n})\}\geq  \frac{\min\{1,r\}}{2^{i_0}}.
	\end{eqnarray*}  Therefore, $(A_n)$ cannot converge to $A$ with respect to  $d_{\mathcal{S}}^{\mathcal{A}_1}$. 
\end{proof}
\begin{remark}\label{metricdiffbornologiesremark}
Let $\tau(d_{\mathcal{S}}^{\mathcal{A}_1})$ and $\tau(d_{\mathcal{C}}^{\mathcal{A}_2})$ be the topologies on $CL(X)$ corresponding to the metrics $d_{\mathcal{S}}^{\mathcal{A}_1}$ and $d_{\mathcal{C}}^{\mathcal{A}_2}$. Then the proof of Theorem \ref{metricdiffbornologies} can be modified to show that  $\tau(d_{\mathcal{C}}^{\mathcal{A}_2}) \subseteq \tau(d_{\mathcal{S}}^{\mathcal{A}_1})$ if and only if $\mathcal{A}_2^* \subseteq \mathcal{A}_1^*$.
\end{remark}
\begin{corollary}\label{metricsamebornologies}
Let $(X,d)$ be a metric space and let $\mathcal{S}$ be a bornology on $X$. Suppose $\mathcal{S}$ is a totally bounded family on $X$, and $\mathcal{A}_1$, $\mathcal{A}_2$ are two countable subcollections of $\mathcal{S}$  satisfying Definition \ref{totally bounded family}. Then the metrics $d_{\mathcal{S}}^{\mathcal{A}_1}$ and $d_{\mathcal{S}}^{\mathcal{A}_2}$ are equivalent on $CL(X)$.
\end{corollary}
\begin{proof}
	One can easily verify that $\mathcal{A}_1^* = \mathcal{A}_2^* = \mathcal{S}^*$. Consequently, by Theorem \ref{metricdiffbornologies}, the metrics $d_{\mathcal{S}}^{\mathcal{A}_1}$ and $d_{\mathcal{S}}^{\mathcal{A}_2}$ are equivalent on $CL(X)$. 
\end{proof}
By Proposition \ref{wijsmancoincidence}, we have $\tau_{\mathcal{S},d} = \tau_{W_{d}}$ if and only if $\mathcal{S} \subseteq \mathcal{TB}_d(X)$. It is natural to ask when the topology generated by metric $d_{\mathcal{S}}^{\mathcal{A}}$ coincides with the Wijsman topology $\tau_{W_d}$? The following theorem answers this question.   
\begin{theorem}\label{dsawijscoincidence}
	Let $(X,d)$ be a metric space and let $\mathcal{S}$
	be a bornology on $X$. Suppose $\mathcal{A} = \{S_n: n \in \mathbb{N}\}$ is a countable subcollection of $\mathcal{S}$ such that $\displaystyle{X = \overline{\bigcup_{S_n \in \mathcal{A}}S_n}}$, then $\tau(d_{\mathcal{S}}^{\mathcal{A}}) = \tau_{W_{d}}$ if and only if $\mathcal{A} \subseteq \mathcal{TB}_d(X)$.  \end{theorem}
\begin{proof}
	Suppose $\tau(d_{\mathcal{S}}^{\mathcal{A}}) = \tau_{W_{d}}$ and $\mathcal{A} \nsubseteq \mathcal{TB}_d(X)$. Then there exists $n_0 \in \mathbb{N}$ such that $S_{n_0} \notin \mathcal{TB}_d(X)$. So we can find an $0 < \epsilon < 1$ such that $S_{n_0} \nsubseteq B_d(F, \epsilon)$ for any $F \in \mathcal{F}(X)$. Direct $\mathcal{F}(X)$ by set inclusion $\subseteq$. Define $A_F = F$ for any $F \in \mathcal{F}(X)$. Then $(A_F)$ is a net in $CL(X)$. It is easy to see that $(A_F)$ is $\tau_{W_{d}}$-convergent to $X$. However,  $(A_F)$ does not converge to $X$ in $(CL(X), \tau(d_{\mathcal{S}}^{\mathcal{A}})) $. To see this, for any $F \in \mathcal{F}(X)$, consider $$d_{\mathcal{S}}^{\mathcal{A}}(A_F, X) = \sum_{n=1}^{\infty}2^{-n}\min\{1, \sup_{x \in S_n} \vert d(x, X) - d(x, A_F)\vert\}\geq \frac{\epsilon}{2^{n_0}}.$$ Hence, we arrive at a contradiction. 
	
	Conversely, suppose $\mathcal{A} \subseteq \mathcal{TB}_d(X)$. Using idea of the proof of Proposition \ref{TBsepimpliessep}, $X$ is separable. Let $\mathcal{A}' = \{\{x_n\}: n \in \mathbb{N}\} \subseteq \mathcal{F}(X)$ such that $X = \overline{\cup \mathcal{A}'}$. Then $\mathcal{A}^* \subseteq \mathcal{A}'^*$. So by Remark \ref{metricdiffbornologiesremark}, $\tau(d_{\mathcal{S}}^{\mathcal{A}}) \subseteq \tau(d_{\mathcal{F}(X)}^{\mathcal{A}'})$. Also by Exercise $2.1.6$ of (\cite{ToCCoS}), $\tau(d_{\mathcal{F}(X)}^{\mathcal{A}'}) = \tau_{W_{d}}$. Therefore, by Proposition \ref{dSAfinerwijsman}, $\tau(d_{\mathcal{S}}^{\mathcal{A}}) = \tau_{W_{d}}$ on $CL(X)$. 
\end{proof}

The following example shows that in general, convergence with respect to metric $d_{\mathcal{S}}^{\mathcal{A}}$ corresponding to a subcollection $\mathcal{A}$ of a bornology $\mathcal{S}$ satisfying $\displaystyle{\overline{\bigcup_{A \in \mathcal{A}}A} = X}$ may not imply $\mathcal{S}$-convergence.  
\begin{example}\label{metricSconvergence}
	Consider $(X,d)$ and $\mathcal{S}$ as in Example \ref{counterexample1}. Let $\mathcal{A}= \{[-n,n]: n \in \mathbb{N}\}$. Then using Proposition \ref{metric}, $d_{\mathcal{S}}^{\mathcal{A}}$ forms a metric on $CL(X)$. For each $n \in \mathbb{N}$, define $A_n = [-n,n]$ and $A = \mathbb{R}$. We show that the sequence $(A_n)$ is $d_{\mathcal{S}}^{\mathcal{A}}$-convergent to $A$ while its $\mathcal{S}$-convergence to $A$ fails. 
	
	Let $\epsilon > 0$. Choose $n_0 \in \mathbb{N}$ such that $\frac{1}{2^{n_0}} < \epsilon$. Consider 
	\begin{eqnarray*}
		d_{\mathcal{S}}^{\mathcal{A}}(A, A_n) = \sum_{m =1 }^{\infty}2^{-m}\min\{1, \sup_{x \in [-m,m]}\vert d(x, A) - d(x, A_n)\vert\} \\ = \sum_{m =1}^{\infty}2^{-m}\min \{1, \sup_{x \in [-m,m]}\vert d(x, A_n)\vert\} \\
		= \sum_{m > n}^{\infty}2^{-m} \min \{1, \sup_{x \in [-m,m]}\vert d(x, A_n)\vert\} = \frac{1}{2^n}.
	\end{eqnarray*}
	So for $n \geq n_0$, we have $d_{\mathcal{S}}^{\mathcal{A}}(A, A_n) \leq \frac{1}{2^{n_0}} <  \epsilon$. Thus, $(A_n)$ converges to $A$ with respect to $d_{\mathcal{S}}^{\mathcal{A}}$.   
	However, for $S = \{n\sqrt{2}: n \in \mathbb{N}\} \in \mathcal{C}(\mathbb{Q}^c)$, $A \cap S \nsubseteq B_d(A_n, \epsilon)$ for any $n \in \mathbb{N}$ and $0< \epsilon < 1$. Thus, the $\mathcal{S}$-convergence of sequence $(A_n)$ to $A$ fails.\qed   \end{example} 
Additionally, it can also be concluded from Example \ref{metricSconvergence} that for a bornology $\mathcal{S}$ on $(X,d)$, the topology generated by the metric $d_{\mathcal{S}}^{\mathcal{A}}$ may not be equal to $\tau_{\mathcal{S}, d}$ as $\tau_{\mathcal{S}, d}$-convergence is always finer than the $\mathcal{S}$-convergence. 
%

 The next result shows that for a bornology $\mathcal{S}$ on a metric space $(X,d)$ such that $X$ is $\mathcal{S}$-separable, the topology corresponding to the metric $d_{\mathcal{S}}^{\mathcal{A}}$ is always weaker than the topology $\tau_{\mathcal{S}, d}$. 
A topological space $(X, \tau)$ is said to have a \textit{$G_{\delta}$-diagonal} if the diagonal $\Delta = \{(x,x ) :  x \in X\}$ is a $G_{\delta}$-set in $X \times X$. It is easy to see that every submetrizable topological space $X$ has a $G_{\delta}$-diagonal, and if $X$ has a $G_{\delta}$-diagonal, then points of $X$ are $G_{\delta}$.  
\begin{theorem}\label{Submetrizability}
	Let $(X,d)$ be a metric space and let $\mathcal{S}$ be a bornology on $X$. Then the following statements are equivalent:
	\begin{enumerate}
		\item[(i)] $X$ is $\mathcal{S}$-separable;
		\item[(ii)] $(CL(X),\tau_{\mathcal{S},d})$ is submetrizable;
		\item[(iii)] $(CL(X),\tau_{\mathcal{S},d})$ has a $G_{\delta}$-diagonal;
		\item[(iv)] each point of $(CL(X),\tau_{\mathcal{S},d})$ is a $G_{\delta}$-set.
		
	\end{enumerate}
\end{theorem}
\begin{proof} It is enough to prove $(i) \Rightarrow (ii)$ and $(iv) \Rightarrow (i)$.\\
	$(i) \Rightarrow (ii)$.
	Suppose $X$ is $\mathcal{S}$-separable. Let $\mathcal{A} = \{S_n: n \in \mathbb{N}\} \subseteq \mathcal{S}$ such that $\displaystyle{X = \overline{\bigcup_{S_n \in \mathcal{A}}S_n}}$. So by Proposition \ref{metric}, $$d_{\mathcal{S}}^{\mathcal{A}}(A,C) = \sum_{n=1}^{\infty}2^{-n}\min\{1, \sup_{x \in S_{n}}\vert d(x,A) - d(x,C) \vert\}$$
	is a metric on $CL(X)$. We show that the metric $d_{\mathcal{S}}^{\mathcal{A}}$ gives a coarser topology than $\tau_{\mathcal{S},d}$ on $CL(X)$. Consider any open ball $B_{d_{\mathcal{S}}^{\mathcal{A}}}(A, \epsilon)$ in $(CL(X), d_{\mathcal{S}}^{\mathcal{A}})$. We show that there exist $S \in \mathcal{S}$ and $r > 0$ such that $[S,r](A) \subseteq B_{d_{\mathcal{S}}^{\mathcal{A}}}(A, \epsilon)$. It can be easily verified that for $\epsilon > 0$, there exists $0 < r < 1$ and $n_{0} \in \mathbb{N}$ such that $$r\left(1- \frac{1}{2^{n_0}}\right) + \frac{1}{2^{n_0}} < \epsilon.$$ Put $\displaystyle{S = \bigcup_{n=1}^{n_0}S_n}.$  Let $C \in [S,r](A)$. Then for each $n\leq n_0$, we have $|d(x,A) - d(x,C)| < r$ for all $x \in S_n$.  Consequently, \begin{eqnarray*}\mathclap{\begin{split}		
	\sum_{n=1}^{\infty} 2^{-n} \min \{1, \sup_{x \in S_{n}}|d(x, A) - d(x, C)|\} = \sum_{n=1}^{n_0} 2^{-n} \min \{1, \sup_{x \in S_{n}}|d(x, A) - d(x, C)|\}+\\ \sum_{n=n_0+1}^{\infty} 2^{-n} \min \{1, \sup_{x \in S_{n}}|d(x, A) - d(x, C)|\}\\\leq \sum_{n=1}^{n_0} 2^{-n} \min \{1, \sup_{x \in S_{n}}|d(x, A) - d(x, C)|\} + \frac{1}{2^{n_0}} < \epsilon.\end{split}}\end{eqnarray*} Hence $C \in B_{d_{\mathcal{S}}^{\mathcal{A}}}(A, \epsilon)$.

	$(iv) \Rightarrow (i)$. Let $\{[S_n,\epsilon_n](X) \vert\ \ n \in \mathbb{N}\}$ be a family of open neighborhoods of $X$ in $(CL(X),\tau_{\mathcal{S},d})$ such that $\displaystyle{\{X\} = \bigcap_{n \in \mathbb{N}}[S_n,\epsilon_n](X)}$. We claim that $\displaystyle{X = \overline{\bigcup_{n \in \mathbb{N}}S_n}}$. Suppose on the contrary, there is a $y \in X$ such that $y \notin \displaystyle{\overline{\bigcup_{n \in \mathbb{N}}S_n}}$. So there exists $\delta > 0$ such that $\displaystyle{B_d\left(y, \frac{\delta}{2}\right) \cap \left(\bigcup_{n \in \mathbb{N}}S_n\right) = \emptyset}$. Set $A = X \setminus B_d(y, \frac{\delta}{2})$. Then $A \in CL(X)$. Now for any $x \in S_n$, we have $d(x,A)= 0$. Thus, $\displaystyle{A \in \bigcap_{n \in \mathbb{N}}[S_n,\epsilon_n](X)}$. We arrive at a contradiction.
	\end{proof}
It is known that the $\mathcal{K}(X)$-convergence is compatible with the classical Fell topology on $CL(X)$ (\cite{beer1987fell}). Beer et al. in (\cite{Bcas}) have shown that whenever $\mathcal{S}$ is shielded from closed sets, the topology of $\mathcal{S}$-convergence is metrizable if and only if $\mathcal{S}$ has a countable base.
 We now give a weaker result about the metrizability of the topology of $\mathcal{S}$-convergence that follows immediately from Theorem \ref{Submetrizability}.
\begin{corollary}
Let $(X,d)$ be a metric space and let $\mathcal{S}$ be a bornology on $X$. If the topology of $\mathcal{S}$-convergence on $CL(X)$ is metrizable, then $X$ is $\mathcal{S}$-separable. In particular, if the Fell topology on $CL(X)$ is metrizable, then $X$ is separable. 
\end{corollary}
\begin{remark}
	It can be inferred from Theorem \ref{Submetrizability} that for a non-separable metric space, there does not exist any metrizable topology weaker than $\tau_{W_{d}}$ on $CL(X)$.
\end{remark}


\section{Metrizability}

 In this section, we explore the conditions for which the space $(CL(X), \tau_{\mathcal{S},d})$ is metrizable or first countable. Since the metrizability of a topological space implies submetrizability, the condition of $X$ being $\mathcal{S}$-separable is necessary for the metrizability of the space $(CL(X), \tau_{\mathcal{S},d})$. However, it may not be sufficient.

%


\begin{example}
	Consider the metric space $(X,d)$ and bornology $\mathcal{S}$ as in Example \ref{counterexample1}. Then $X$ is $\mathcal{S}$-separable. We show that the space $(CL(X), \tau_{\mathcal{S},d})$ is not metrizable. Suppose $(CL(X), \tau_{\mathcal{S},d})$ is metrizable. Then each $C \in CL(X)$ has a countable local base in $(CL(X), \tau_{\mathcal{S},d})$. In particular, let $\{[S_n, \frac{1}{m}](X): S_n \in \mathcal{S}, n,m \in \mathbb{N}\}$ be a countable local base at $X$. We can assume that each $S_n$ is of the form $[-n,n]\cup C_n$, where $C_n \in \mathcal{C}(\mathbb{Q}^c)$. Let $\displaystyle{D = \mathbb{Q}^c \setminus \left(\bigcup_{n \in \mathbb{N}}C_n\right)}$. Then by proceeding as in Example \ref{counterexample1}, we can find $B \in \mathcal{S}$ such that $B = \{a_n: n \in \mathbb{N}\}\cup\{b_n:  n \in \mathbb{N}\}$, where for each $n \in \mathbb{N}$, $a_n < -n$, $b_n > n$ and $B \subseteq D$.  Take $0 < \epsilon < 1$. Then $[B, \epsilon](X)$ is a $\tau_{\mathcal{S},d}$-neighborhood of $X$. However, for any $n,m \in \mathbb{N}$, $[S_n, \frac{1}{m}](X)$ is not contained in $[B, \epsilon](X)$. Since for any  $n \in \mathbb{N}$, $S_n= \overline{S_n} \in [S_n, \frac{1}{m}](X)$ and $S_n \notin [B, \epsilon](X)$. Which is a contradiction.\qed
\end{example}

We next give a sufficient condition for the metrizability of the hyperspace, $(CL(X), \tau_{\mathcal{S},d})$. In fact, under the same condition the compatible uniformity $\mathcal{U}_{\mathcal{S},d}$ for $(CL(X), \tau_{\mathcal{S},d})$ is also metrizable. 

\begin{theorem}\label{metrizableuniformity}
	Let $(X,d)$ be a metric space and let $\mathcal{S}$ be a bornology on $X$. If $\mathcal{S}$ is a totally bounded family on $X$, then $(CL(X),~ \mathcal{U}_{\mathcal{S},d})$ is metrizable.
\end{theorem}

\begin{proof}
	We show that the uniformity $\mathcal{U}_{\mathcal{S},d}$ has a countable base. Given $\mathcal{S}$ is a totally bounded family on $X$, let $\{S_n: n \in \mathbb{N}\} \subseteq \mathcal{S}$  be such that $S_n \subseteq S_{n+1}$ for each $n$ and that it satisfies Definition \ref{totally bounded family}. We claim that $\{[S_n, \frac{1}{m}]: n,m \in \mathbb{N}\}$ forms a countable base for the uniformity $\mathcal{U}_{\mathcal{S},d}$. Pick $S \in \mathcal{S}$ and $\epsilon > 0$. Then there is an $n_0 \in \mathbb{N}$ such that $\frac{1}{n_0} < \frac{\epsilon}{8}$ and $S \subseteq B_d(S_{n_{0}}, \frac{\epsilon}{8})$. Let $x \in S$. Consequently, $d(x,y) < \frac{\epsilon}{8}$ for some $y \in S_{n_{0}}$. Consider $(A,C) \in [S_{n_{0}}, \frac{1}{n_0}]$. Then $$d(x,A) \leq d(x,y) + d(y, A) \leq d(y, C) + \frac{\epsilon}{4} < d(x,C) + \frac{\epsilon}{2}. $$
	Similarly, $$d(x,C) \leq d(x,y) + d(y,C) \leq d(y, A) + \frac{\epsilon}{4} < d(x,A) + \frac{\epsilon}{2}. $$
  So $\displaystyle{\sup_{x \in S}\vert d(x, A)- d(x, C)\vert < \epsilon}$, that is, $(A,C) \in [S, \epsilon]$. Thus, $[S_{n_{0}}, \frac{1}{n_0}] \subseteq [S, \epsilon]$. 
	\end{proof}
\begin{corollary}
	Let $(X,d)$ be a metric space and let $\mathcal{S}$ be a bornology on $X$. If $\mathcal{S}$ is a totally bounded family on $X$, then $(CL(X), \tau_{\mathcal{S},d})$ is metrizable. 
\end{corollary}

We now give a compatible metric for the hyperspace $(CL(X), \tau_{\mathcal{S},d})$ whenever $\mathcal{S}$ is a totally bounded family on $X$. Recall that if $X$ is $\mathcal{S}$-separable, then the metric $d_{\mathcal{S}}^{\mathcal{A}}$ on $CL(X)$ with respect to a countable subcollection $\mathcal{A}$ of $\mathcal{S}$ satisfying Definition \ref{Sseperable} gives a coarser topology than $\tau_{\mathcal{S},d}$ (Theorem \ref{Submetrizability}). The following theorem shows that the metric $d_{\mathcal{S}}^{\mathcal{A}}$ is compatible with the topology $\tau_{\mathcal{S},d}$ whenever $\mathcal{S}$ is a totally bounded family on $X$. 

\begin{theorem}\label{compatiblemetric}
	Let $(X,d)$ be a metric space and let $\mathcal{S}$ be a bornology on $X$. Suppose $\mathcal{S}$ is a totally bounded family on $X$ and $\mathcal{A}$ is any countable subcollection of $\mathcal{S}$ satisfying Definition \ref{totally bounded family}. Then the metric $d_{\mathcal{S}}^{\mathcal{A}}$ is a compatible metric for the hyperspace $(CL(X), \tau_{\mathcal{S},d})$. 
\end{theorem}

\begin{proof}
	Let $\mathcal{A} = \{S_n: n \in \mathbb{N}\} \subseteq \mathcal{S}$ satisfy Definition \ref{totally bounded family}. By Proposition \ref{dSAequivalentmetric}, we can assume $S_n \subseteq S_{n+1}$ for each $n \in \mathbb{N}$. By Theorem \ref{Submetrizability}, the topology generated by the metric $d_{\mathcal{S}}^{\mathcal{A}}$ is weaker than $\tau_{\mathcal{S},d}$ on $CL(X)$. For the reverse inclusion, pick $S \in \mathcal{S}$, $A \in CL(X)$, and $0 < \epsilon < 1$. Choose $N \in \mathbb{N}$ such that $S \subseteq B_d(S_N, \frac{\epsilon}{4})$.  Set $r = \frac{\epsilon}{2^{N+3}}$. We claim that $B_{d_{\mathcal{S}}^{\mathcal{A}}}(A,r) \subseteq [S, \epsilon](A)$. Consider $C \in B_{d_{\mathcal{S}}^{\mathcal{A}}}(A, r)$. Then 
	$$\sum_{n = 1}^{\infty}2^{-n}\min\{1, \sup_{x \in S_n}\vert d(x,A) - d(x, C) \vert \} < r.$$      
	Consequently, by the choice of $r$, we have 
	$$\sup_{x \in S_N}\vert d(x,A) - d(x,C) \vert < \frac{\epsilon}{4}.$$  
	If $x \in S$, then there is a $y \in S_N$ such that $d(x,y) < \frac{\epsilon}{4}$. Then $d(x, A) \leq d(x,y) + d(y,C)+ \frac{\epsilon}{4} \leq d(x,C) + \frac{3\epsilon}{4}$. Similarly, $d(x,C) \leq d(x,y) + d(y,C) \leq d(x,A) + \frac{3\epsilon}{4}$. Therefore, 
	$C \in [S, \epsilon] (A)$. 
\end{proof}
Since every bornology $\mathcal{S}$ with a countable base is a totally bounded family on $X$, Theorem \ref{metrizableuniformity} gives a weaker condition for the metrizability of $(CL(X), \tau_{\mathcal{S},d})$ than Proposition 3 of (\cite{Idealtopo}). Moreover, Theorem \ref{compatiblemetric} provides a compatible metric for $(CL(X),\tau_{\mathcal{S},d})$ under the similar assumption.

The following result is an immediate corollary to Theorem \ref{metricdiffbornologies} and Theorem \ref{compatiblemetric}.

	\begin{corollary}\label{comptiblemetric_corollary 1}
	Let $(X,d)$ be a metric space and let $\mathcal{S}$, $\mathcal{C}$ be two bornologies on $X$. Suppose $\mathcal{A}_1$ and $\mathcal{A}_2$ are countable subcollections of $\mathcal{S}$ and $\mathcal{C}$, respectively satisfying Definition \ref{totally bounded family}. Then the following statements are equivalent:
	\begin{enumerate}[(i)]
		\item metrics $d_{\mathcal{S}}^{\mathcal{A}_1}$ and $d_{\mathcal{C}}^{\mathcal{A}_2}$ are equivalent; 
		\item $\tau_{\mathcal{S}, d} = \tau_{\mathcal{C},d}$;
		\item $\mathcal{A}_1^* = \mathcal{A}_2^*$. 
	\end{enumerate} 
\end{corollary}

\begin{remark}

	 It is to be noted that for two different countable subcollections $\mathcal{A}_1$ and $\mathcal{A}_2$ satisfying Definition \ref{Sseperable}, the topologies generated by the metrics $d_{\mathcal{S}}^{\mathcal{A}_1}$ and $d_{\mathcal{S}}^{\mathcal{A}_2}$ need not be equivalent (see, Example \ref{metric_diffcollection}). However, for a totally bounded family $\mathcal{S}$ on $X$ whenever $\mathcal{A}_1$ and $\mathcal{A}_2$ satisfy Definition \ref{totally bounded family}, by Theorem \ref{compatiblemetric}, the topologies generated by the metrics $d_{\mathcal{S}}^{\mathcal{A}_1}$ and $d_{\mathcal{S}}^{\mathcal{A}_2}$ are equal to $\tau_{\mathcal{S}, d}$. Consequently, in this case the convergence given by the metric $d_{\mathcal{S}}^{\mathcal{A}}$ is finer than  the $\mathcal{S}$-convergence.
\end{remark}


\begin{corollary} \normalfont{(\cite{ToCCoS})}\label{comptiblemetric_corollary 2}
	Let $(X,d)$ be a metric space and let $D = \{x_n: n \in \mathbb{N}\}$ be a countable dense subset of $X$. Then there exists a 
	compatible metric for the space $(CL(X), \tau_{W_{d}})$.   
\end{corollary} 

\begin{proof}
	Let $\mathcal{A} = \{\{x_n\}: n \in \mathbb{N}\}$. Then $\mathcal{A}$ satisfies Definition \ref{totally bounded family} for $\mathcal{F}(X)$. Consequently, by Theorem \ref{compatiblemetric}, the metric $d_{\mathcal{F}(X)}^{\mathcal{A}}$ defined by $$d_{\mathcal{F}(X)}^{\mathcal{A}}(A,C) = \sum_{n = 1}^{\infty}2^{-n}\min\{1, \vert d(x_n, A) - d(x_n, C) \vert\},$$ where $A,C \in CL(X)$, is the required compatible metric for $(CL(X), \tau_{W_{d}})$. \end{proof}

\begin{corollary}\normalfont{(\cite{ToCCoS})}\label{comptiblemetric_corollary 3}
	Let $(X,d)$ be a metric space. Then there exists a
	 compatible metric for $(CL(X), \tau_{AW_{d}})$.   
\end{corollary} 

\begin{proof}
Let $x_0 \in X$. Then $\mathcal{B}_d(X)$ has a countable base $\mathcal{A} = \{B_d(x_0,n): n \in \mathbb{N}\}$. Consequently, by Theorem \ref{compatiblemetric}, the metric $d_{\mathcal{B}_d(X)}^{\mathcal{A}}$ defined by $$d_{\mathcal{B}_d(X)}^{\mathcal{A}}(A,C) = \sum_{n = 1}^{\infty}2^{-n}\min\{1, \sup_{d(x,x_0) < n} \vert d(x, A) - d(x, C) \vert\},$$ where $A,C \in CL(X)$, is the required compatible metric for $(CL(X), \tau_{AW_{d}})$.\end{proof}

 Finally, for $\mathcal{S} = \mathcal{P}_0(X)$ and $\mathcal{A}=\{X\}$, we have $d_{\mathcal{P}_0(X)}^{\mathcal{A}}(A,C) = \frac{1}{2}\min\{1, H_d(A,C)\}$, where $A,C \in CL(X)$. Therefore, the metric $d_{\mathcal{P}_0(X)}^{\mathcal{A}}$ is a compatible metric for the Hausdorff distance topology on $CL(X)$. 

 Our next corollary can be deduced using Corollaries \ref{comptiblemetric_corollary 1}-\ref{comptiblemetric_corollary 3} and the above discussion. One can easily verify the following facts: 
 \begin{itemize}
 	\item if $D = \{x_n : n \in \mathbb{N}\}$
 	is dense in $(X,d)$ and $\mathcal{A}_1 = \{\{x_n\}: n \in \mathbb{N}\}$,  then $\mathcal{A}_1^* = \mathcal{TB}_d(X)$;
 	\item if $x_0 \in X$ and $\mathcal{A}_2 = \{B_d(x_0, n):  n \in \mathbb{N}\}$, then $\mathcal{A}_2^*  =\mathcal{B}_d(X)$;
 	\item if $\mathcal{A}_3 = \{X\}$, then $\mathcal{A}_3^* = \mathcal{P}_0(X)$.    
 \end{itemize}
 
 \begin{corollary}\normalfont{(}\cite{ToCCoS}\normalfont{)}
 	Let $(X,d)$ be a metric space. Then the following statements are equivalent:
 	\begin{enumerate}[(i)]
 		\item $\tau_{AW_{d}} = \tau_{W_{d}}$ on $CL(X)$ if and only if $\mathcal{B}_d(X) = \mathcal{TB}_d(X)$;
 		\item $\tau_{H_{d}} = \tau_{W_{d}}$ on $CL(X)$ if and only if $\mathcal{P}_0(X) = \mathcal{TB}_d(X)$;
 		\item $\tau_{H_{d}} =\tau_{AW_{d}}$ on $CL(X)$ if and only if $\mathcal{P}_0(X) = \mathcal{B}_d(X)$.
 	\end{enumerate} 
 \end{corollary}
In the following result, we show that the condition of $\mathcal{S}$ being a totally bounded family on $X$ is also necessary for the metrizability of the space $(CL(X), \tau_{\mathcal{S},d})$.

\begin{theorem}\label{metrizabilityiff}
	Let $(X,d)$ be a metric space and let $\mathcal{S}$ be a bornology on $X$. Then the following statements are equivalent:
	\begin{enumerate}[(i)]
		\item $\mathcal{S}$ is a totally bounded family on  $X$;
		
		\item $(CL(X), \mathcal{U}_{\mathcal{S},d})$ is metrizable; 
		
		\item $(CL(X), \tau_{\mathcal{S},d})$ is metrizable;  
		\item  $(CL(X), \tau_{\mathcal{S},d})$ is first countable.
	\end{enumerate}
\end{theorem}  
\begin{proof}
	$(i) \Rightarrow (ii)$. It follows from Theorem \ref{metrizableuniformity}. Also $(ii) \Rightarrow (iii) \Rightarrow (iv)$ are immediate. 
	
	$(iv) \Rightarrow (i)$. Suppose by contradiction that $(i)$ fails. By the hypothesis, $X$ has a countable local base in $(CL(X), \tau_{\mathcal{S},d})$. Without loss of generality, assume $\{[S_n, \frac{1}{m}](X): n,m \in \mathbb{N}\}$ is a countable local base at $X$ in $(CL(X), \tau_{\mathcal{S}, d})$, where $S_n \subseteq S_{n+1}$. Since $\mathcal{S}$ is not a totally bounded family on $X$, for the countable subcollection $\{S_n: n \in \mathbb{N}\}$ of $\mathcal{S}$, there is an $S \in \mathcal{S}$ and $\epsilon > 0$ such that $S \nsubseteq B_d(S_n, \epsilon)$ for any $n \in \mathbb{N}$. Consider the neighborhood $\displaystyle{[S, \epsilon](X) = \{C \in CL(X): \sup_{x\in S}d(x, C) < \epsilon\}}$ of $X$ in $(CL(X), \tau_{\mathcal{S},d})$. Observe that for each $n,m \in \mathbb{N}$, $\overline{S_n} \in [S_n, \frac{1}{m}](X)$ while $\overline{S_n} \notin [S, \epsilon](X)$. We arrive at a contradiction.   
	\end{proof}

\begin{corollary}\normalfont{(\cite{ToCCoS})}
	Let $(X,d)$ be a metric space. Then the following statements are equivalent:
	\begin{enumerate}[(i)]
		\item $(X,d)$ is separable;
		\item $(CL(X), \tau_{W_{d}})$ is metrizable;
		\item $(CL(X), \tau_{W_d})$ is first countable.
		\end{enumerate}
	\end{corollary}



\section{Countability Properties}

In this section, we investigate separability and second countability  of the hyperspace $(CL(X), \tau_{\mathcal{S},d})$.


\begin{proposition}\label{2ndcountabilitynecessary}
	Let $(X,d)$ be a metric space and let $\mathcal{S}$ be a bornology on $X$. If $(CL(X), \tau_{\mathcal{S},d})$ is second countable, then $\mathcal{S}$ is a totally bounded family on $X$. 
\end{proposition}

\begin{proof}
	Since every second countable space is first countable. So by Theorem \ref{metrizabilityiff}, the result follows.  
\end{proof}

The previous result tells that $\mathcal{S}$ being a totally bounded family on $X$ is a necessary condition for the second countability of the hyperspace $(CL(X), \tau_{\mathcal{S},d})$. However, this may not be sufficient for second countability of $(CL(X), \tau_{\mathcal{S},d})$. For example, $\mathcal{B}_d(X)$ is always a totally bounded family on $X$ but $(CL(X), \tau_{AW_{d}})$ need not be second countable unless $\mathcal{B}_d(X) = \mathcal{TB}_d(X)$. 

\begin{proposition}\label{separabilityTsd}
	Let $(X,d)$ be a metric space  and let $\mathcal{S}$ be a bornology on $X$. If $(CL(X), \tau_{\mathcal{S},d})$ is separable, then $\mathcal{S} \subseteq \mathcal{TB}_d(X)$. 
\end{proposition}
\begin{proof}
	Suppose $\mathcal{S} \nsubseteq \mathcal{T}\mathcal{B}_d(X)$. Choose $S \in \mathcal{S}$ such that $S$ is not totally bounded. Then for some $\epsilon > 0$, we can find an infinite $\epsilon$-discrete subset $E$ of $S$, that is, $d(x,y) > \epsilon$ for any $x,y \in E$. Let $\mathscr{A} = \{D: D \subset E\}$. We claim that for $0 < r < \frac{\epsilon}{4}$, the family $\{[S,r](D): D \in \mathscr{A}\}$ is an uncountable pairwise disjoint collection of open sets in $(CL(X), \tau_{\mathcal{S},d})$. Let $D_1, D_2$ be two distinct elements of $\mathscr{A}$. By contradiction, suppose $A \in [S,r](D_1) \cap [S,r](D_2)$. For $x \in D_1\setminus D_2$, we have $d(x,A) < r$. Since $A \in [S,r](D_2)$, we get $d(x,D_2) < 2r$. However, by the $\epsilon$-discreteness of $E$, we have $d(x, D_2) > 4r$. We arrive at a contradiction. Therefore, $(CL(X), \tau_{\mathcal{S},d})$ is not separable.   
\end{proof}
 \begin{theorem}
 	Let $(X,d)$ be a metric space and let $\mathcal{S}$ be a bornology on $X$. If $X$ is $\mathcal{S}$-separable, then the following statements are equivalent:
 	\begin{enumerate}[(i)]
 		\item $\mathcal{S} \subseteq \mathcal{TB}_d(X)$;
 		\item $(CL(X), \tau_{\mathcal{S}, d})$ is separable.
\end{enumerate} 
\end{theorem}
\begin{proof}
	$(i)\Rightarrow (ii)$.  Given $\mathcal{S} \subseteq \mathcal{T}\mathcal{B}_d(X)$, by Proposition \ref{wijsmancoincidence}, we have $\tau_{\mathcal{S},d} = \tau_{W_{d}}$. Since $\mathcal{S} \subseteq \mathcal{T}\mathcal{B}_d(X)$ and $X$ is $\mathcal{S}$-separable, by Proposition \ref{TBsepimpliessep}, $X$ is separable. So by Theorem $2.1.5$ of (\cite{ToCCoS}), $(CL(X),\tau_{W_{d}})$ is second countable. Thus, $(ii)$ follows. 
	
	$(ii) \Rightarrow (i)$. It follows from  Proposition \ref{separabilityTsd}. 
\end{proof}
\begin{theorem}
	Let $(X,d)$ be a metric space and let $\mathcal{S}$ be a bornology on $X$. If $\mathcal{A} = \{S_n:n \in \mathbb{N}\}$ is a countable subcollection of $\mathcal{S}$  such that $\displaystyle{X = \overline{\bigcup_{S_n \in \mathcal{A}}S_n}}$. Then the following statements are equivalent:
	\begin{enumerate}[(i)]
		\item $\mathcal{A} \subseteq \mathcal{TB}_d(X)$;
		\item $(CL(X), d_{\mathcal{S}}^{\mathcal{A}})$ is separable.
	\end{enumerate}
\end{theorem}
\begin{proof}
	
	$(i) \Rightarrow (ii)$. Given $X$ is $\mathcal{S}$-separable and $\mathcal{A} \subseteq \mathcal{TB}_d(X)$, by using the idea of proof of Proposition \ref{TBsepimpliessep}, we can show that $X$ is separable. Consequently, by Theorem $2.1.5$ of (\cite{ToCCoS}), $(CL(X), \tau_{W_d})$ is separable. But by Theorem \ref{dsawijscoincidence}, we have $\tau(d_{\mathcal{S}}^{\mathcal{A}}) = \tau_{W_d}$.  
	
	$(ii) \Rightarrow (i)$. Suppose $\mathcal{A} \nsubseteq \mathcal{T}\mathcal{B}_d(X)$. Choose $S_{n_0} \in \mathcal{A}$ such that $S_{n_0}$ is not totally bounded. Then for some $\epsilon > 0$, we can find an infinite $\epsilon$-discrete subset $E$ of $S_{n_0}$, that is, $d(z_1,z_2) > \epsilon$ for any $z_1,z_2 \in E$. Let $\mathscr{D} = \{D : D\subseteq E\}$.  Let $D, D' \in \mathscr{D}$ be distinct. Suppose $x_0\in D\setminus D'$. Then
	 \begin{eqnarray*}\begin{split}
			d_{\mathcal{S}}^{\mathcal{A}}(D, D')= \sum_{n=1}^{\infty}2^{-n}\min\{1, \sup_{x \in S_{n}}|d(x,D) - d(x,D')| \}
			\\ \geq 2^{-n_{0}}\min\{1, \sup_{x \in S_{n_{0}}}|d(x,D) - d(x,D')|\} \\ \geq 2^{-n_{0}}\min\{1, |d(x_0,D) - d(x_0,D')|\}\\  \geq  2^{-n_{0}}\min\{1, \epsilon\}.\end{split}\end{eqnarray*}
	Thus, $\mathscr{D}$ is an uncountable uniformly discrete subset of $(CL(X), d_{\mathcal{S}}^{\mathcal{A}})$. Hence $(CL(X), d_{\mathcal{S}}^{\mathcal{A}})$  cannot be separable. 
\end{proof}
\begin{theorem}\label{secondcountability}
	Let $(X,d)$ be a metric space and let $\mathcal{S}$ be a bornology on $X$. Then the following statements are equivalent:
	\begin{enumerate}[(i)]
		\item $\mathcal{S}$ is a totally bounded family on $X$ and $\mathcal{S} \subseteq \mathcal{T}\mathcal{B}_d(X)$;
		
		\item $(CL(X), \tau_{\mathcal{S},d})$ is second countable. 
	\end{enumerate}
\end{theorem}

\begin{proof}
$(i) \Rightarrow (ii)$. Given $\mathcal{S} \subseteq \mathcal{T}\mathcal{B}_d(X)$, by Proposition \ref{wijsmancoincidence}, we have $\tau_{\mathcal{S},d} = \tau_{W_{d}}$. Also whenever $\mathcal{S} \subseteq \mathcal{T}\mathcal{B}_d(X)$ and $\mathcal{S}$ is a totally bounded family on $X$, then $X$ is separable. Thus, by Theorem $2.1.5$ of (\cite{ToCCoS}), $(ii)$ is established. 

$(ii) \Rightarrow (i)$. It follows from Propositions \ref{2ndcountabilitynecessary} and \ref{separabilityTsd}.  
\end{proof}

 
 Next we deduce results related to classical topologies as a direct consequence of our Theorem \ref{secondcountability} and Proposition \ref{wijsmancoincidence}. 
\begin{corollary}\normalfont{(\cite{ToCCoS})}
	Let $(X,d)$ be  a metric space. Then the following statements are equivalent: 
	
	\begin{enumerate}[(i)]
		\item $X$ is separable; 
		\item $(CL(X), \tau_{W_{d}})$ is second countable;
	\end{enumerate}
\end{corollary}

\begin{corollary}\normalfont{(\cite{ToCCoS})}
	Let $(X,d)$ be  a metric space. Then the following statements are equivalent:
	\begin{enumerate}[(i)]
		\item closed and bounded sets are totally bounded; 
		\item $\tau_{AW_{d}} =\tau_{W_{d}}$ on $CL(X)$;
		\item $(CL(X), \tau_{AW_{d}})$ is second countable;
		
	\end{enumerate}
\end{corollary}

\begin{corollary}\normalfont{(\cite{ToCCoS})}
	Let $(X,d)$ be  a metric space. Then the following statements are equivalent:
	\begin{enumerate}[(i)]
		\item $X$ is totally bounded; 
		\item $\tau_{H_{d}} =\tau_{W_{d}}$ on $CL(X)$;
		\item $(CL(X), \tau_{H_{d}})$ is second countable.
		
	\end{enumerate}
\end{corollary}
 We now give a sufficient condition for second countability of the topology $\tau(\mathcal{S})$ of $\mathcal{S}$-convergence on $CL(X)$.
\begin{proposition}
	Let $(X,d)$ be a metric space and let $\mathcal{S}$  be a bornology on $X$ such that $X$ is $\mathcal{S}$-separable and $\mathcal{S} \subseteq \mathcal{T}\mathcal{B}_d(X)$. If each proper closed ball of $X$ lies in $\mathcal{S}$, then $(CL(X), \tau(\mathcal{S}))$ is second countable.
\end{proposition}

\begin{proof}
	Given each proper closed ball is in $\mathcal{S}$, so by Theorem $34.3$ of (\cite{beer2023bornologies}), we have $\mathcal{S}$-convergence implies $\tau_{W_{d}}$-convergence.  Since $\mathcal{S} \subseteq \mathcal{T}\mathcal{B}_d(X)$, we have $\tau_{\mathcal{S},d} = \tau_{W_{d}}$ on $CL(X)$.  Hence $\tau(\mathcal{S}) = \tau_{\mathcal{S},d}$ under the given hypothesis. But $X$ being $\mathcal{S}$-separable and $\mathcal{S} \subseteq \mathcal{T}\mathcal{B}_d(X)$ implies $\mathcal{S}$ is a totally bounded family on $X$. Therefore, by Theorem \ref{secondcountability}, $(CL(X), \tau(\mathcal{S}))$ is second countable. 
\end{proof}
\begin{corollary}
	Let $(X,d)$ be a metric space. If $X$ is countable and discrete, then the topology of $\mathcal{F}(X)$-convergence on $CL(X)$ is second countable.
\end{corollary}

\bibliographystyle{plain}
\bibliography{reference_file}
\end{document}